\documentclass[11pt]{article}

\usepackage{latexsym,amssymb}
\usepackage{calc}

\newcounter{hours}\newcounter{minutes}

\def\cl{{\rm Cl \,}}

\def\nr{\par \noindent}

\def\Def{\stackrel{\mathrm{def}}{=}}

\def\inter{{\rm int \,}}

\def\E{\mathbb{E}}
\def\H{\mathbb{H}}
\def\R{\mathbb{R}}
\def\S{\mathbb{S}}
\def\V{\mathbb{V}}
\def\beq{\begin{equation}}
\def\eeq{\end{equation}}

\newcommand{\SetEQ}{\setcounter{equation}{0}}
\newcommand{\refLE}[1]{\ensuremath{\stackrel{(\ref{#1})}{\leq}}}
\newcommand{\refEQ}[1]{\ensuremath{\stackrel{(\ref{#1})}{=}}}
\newcommand{\refGE}[1]{\ensuremath{\stackrel{(\ref{#1})}{\geq}}}

\newtheorem{theorem}{Theorem}
\newtheorem{lemma}{Lemma}
\newtheorem{corollary}{Corollary}

\newtheorem{assumption}{Assumption}
\newtheorem{definition}{Definition}

\newtheorem{example}{Example}
\newtheorem{remark}{Remark}
\newcommand{\proof}{\bf Proof: \rm \nr}
\newcommand{\qed}{\hfill $\Box$ \nr \medskip}
\newcommand{\half}{\mbox{${1 \over 2}$}}

\newcommand{\iprod}[2]{\left\langle {#1}, {#2} \right\rangle}

\def\ba{\begin{array}}
\def\ea{\end{array}}
\def\beann{\begin{eqnarray*}}
\def\eeann{\end{eqnarray*}}
\def\bea{\begin{eqnarray}}
\def\eea{\end{eqnarray}}

\def\BT{\begin{theorem}}
\def\ET{\end{theorem}}
\def\BL{\begin{lemma}}
\def\EL{\end{lemma}}
\def\BC{\begin{corollary}}
\def\EC{\end{corollary}}
\def\BE{\begin{example}}
\def\EE{\end{example}}
\def\BD{\begin{definition}}
\def\ED{\end{definition}}
\def\BR{\begin{remark}}
\def\ER{\end{remark}}
\def\BAS{\begin{assumption}}
\def\EAS{\end{assumption}}
\def\BI{\begin{itemize}}
\def\EI{\end{itemize}}

\def\BMP{\begin{minipage}{9.5cm}}
\def\EMP{\end{minipage}}
\def\MPT{\begin{minipage}{11.5cm}}
\def\EPT{\end{minipage}}

\def\la{\langle}
\def\ra{\rangle}

\title{
\textbf{Local Superlinear
Convergence of\\
Polynomial-Time Interior-Point Methods\\
for Hyperbolic Cone Optimization Problems} }
\author{Yu. Nesterov\thanks{Center for Operations Research and Econometrics (CORE),
Catholic University of Louvain (UCL),
\newline \vspace{1ex}
34 voie du Roman Pays, 1348 Louvain-la-Neuve, Belgium;
e-mail: nesterov@core.ucl.ac.be.
\newline
The research results presented in this paper have been
supported by a grant ``Action de recherche concert\`e
\newline ARC 04/09-315'' from the ``Direction de la recherche
scientifique - Communaut\`e fran\c{c}aise de Belgique''.
\newline
The scientific responsibility rests with its author(s).}
 $\mbox{ and}$ L. Tun\c{c}el\thanks{Department of
Combinatorics and Optimization, Faculty
of Mathematics, University of Waterloo, Waterloo, Ontario N2L 3G1,
Canada; e-mail: ltuncel@uwaterloo.ca. Research of this author
was supported in part by Discovery Grants from NSERC.}}

\date{November 2009,
revised: December 2014
}
\begin{document}
\maketitle

\abstract{In this paper, we establish the local superlinear
convergence property of some polynomial-time interior-point methods
for an important family of conic optimization problems.
The main structural property used in our analysis is the
logarithmic homogeneity of self-concordant barrier
function, which must have {\em negative curvature}. We
propose a new path-following predictor-corrector scheme,
which work only in the dual space. It is based on an
easily computable gradient proximity measure, which
ensures an automatic transformation of the global linear
rate of convergence to the local superlinear one under
some mild assumptions. Our step-size procedure for the
predictor step is related to the maximum step size
maintaining feasibility. As the optimal solution set is
approached, our algorithm automatically tightens the
neighborhood of the central path proportionally to the
current duality gap.}

\vspace{3ex}\noindent

\newpage

\section{Introduction}\label{sc-Intro}
\setcounter{equation}{0}

\vspace{1ex}\noindent
{\bf Motivation.} Local superlinear convergence is a
natural and very desirable property of many methods in
Nonlinear Optimization. However, for interior-point
methods the corresponding analysis is not trivial. The
reason is that the barrier function is not defined in a
neighborhood of the solution. Therefore, in order to study
the behavior of the central path, we need to employ
somehow the separable structure of the functional
inequality constraints. From the very beginning \cite{FM},
this analysis was based on the Implicit Function Theorem
as applied to Karush-Kuhn-Tucker conditions.

This tradition explains, to some extend, the delay in
developing an appropriate framework for analyzing the
local behavior of general polynomial-time interior-point
methods \cite{NN}. Indeed, in the theory of
self-concordant functions it is difficult to analyze the
local structure of the solution since we have no access to
the components of the barrier function. Moreover, in
general, it is difficult to relate the self-concordant
barrier to the particular functional inequality
constraints of the initial optimization problem.
Therefore, up to now, the local superlinear convergence
for polynomial-time path-following methods was proved only
for Linear Programming \cite{YGTZ1993,Meh1993} and for
Semidefinite Programming problems
\cite{KSS1998,PS1998,LSZ1998,JPS1999,PS2000}. In both
cases, the authors use in their analysis the special
boundary structure of the feasible regions and of the set
of optimal solutions.

In this paper, we establish the local superlinear convergence
property of interior-point path-following methods by employing some
geometric properties of quite general conic optimization
problem. The main structural property used in our analysis
is the logarithmic homogeneity of self-concordant barrier
functions, and the condition that the barrier must have
negative curvature. We propose a new path-following
predictor-corrector scheme, which works in the dual
space. It is based on an easily computable gradient
proximity measure, which ensures an automatic
transformation of the global linear rate of convergence to
the local superlinear rate (under a mild assumption). Our
step-size procedure for the predictor step is related to
the maximum step size maintaining feasibility. As the
iterates approach the optimal solution set, our algorithm
automatically tightens the neighborhood of the central
path proportionally to the current duality gap.  In the
literature (as we noted above) similar conditions have
been imposed on the iterates of interior-point algorithms
in order to attain local superlinear convergence in the
Semidefinite Programming setting.  However, there exist
different, weaker combinations of conditions that have
been studied in the Semidefinite Programming setting as
well (see for instance
\cite{PS2000,JPS1999,LuMo2007,LuMo2004}).
As a key feature of our analysis, we avoid distinguishing individual
subvectors as ``large'' or ``small.''
Many of the ingredients of our approach are primal-dual symetric;
however, we break the symmetry, when necessary, by our choice of
assumptions and algorithms.  Indeed, in general (beyond the special
case of symmetric cones and self-scaled barriers), only one of the
primal, dual problems admit a logarithmically homogeneous
self-concordant barrier with negative curvature.

\vspace{1ex}\noindent
{\bf Contents.} The paper is organized as follows. In
Section \ref{sc-Pred} we introduce a conic primal-dual
problem and define the central path. After that, we
consider a small full-dimensional dual problem and define
the prediction operator. We derive some representations,
which help to bound the curvature of the central path. In
Section \ref{sc-Meas} we justify the choice of the fixed
Euclidean metrics in the primal and dual spaces. Then, in
Section \ref{sc-MainAss} we introduce two main assumptions
ensuring the quadratic drop in the duality gap for
predicted points from a tight neighborhood of the central
path. The first one is on the strict dual maximum, and the
second one is on the boundedness of the vector $\nabla^2
F_*(s)s_*$ along the central path. The main result of this
section is Theorem \ref{th-DEst} which demonstrates the
quadratic decrease of the distance to the optimal solution
for the prediction point, measured in an appropriately
chosen fixed Euclidean norm.
In Section \ref{sc-PStep}, we estimate efficiency of the
predictor step measured in a local norm defined by the
dual barrier function. Also, we show that the local
quadratic convergence can be achieved by a feasible
predictor step.

In Section \ref{sc-Rec} we prepare for analysis of
polynomial-time predictor-corrector strategies. For that,
we study an important class of barriers with {\em negative
curvature}. This class includes at least the self-scaled
barriers \cite{NT1} and hyperbolic barriers
\cite{Guler1997,Bauschke2001,Renegar2006}, possibly many others.
In Section \ref{sc-Growth} we establish some bounds on the
growth of a variant of the gradient proximity measure. We
show that we can achieve a local superlinear rate of
convergence. It is important to relate the decrease of
penalty parameter of the central path with the
{\em distance to the boundary} of the feasible set while performing the
predictor step. At the same time, we show that, for
achieving the local superlinear convergence, the centering
condition must be satisfied with increasing accuracy.

In Section \ref{sc-Poly} we show that the local
superlinear convergence can be combined with the global
polynomial-time complexity. We present a method, which
works for the barriers with negative curvature, and has a
cheap computation of the predictor step. Finally, in Section
\ref{sc-Disc}, we discuss the results and study two
2D-examples, which demonstrate that our assumptions are
quite natural.

\vspace{1ex}\noindent
{\bf Notation and generalities.} In what follows, we
denote by $\E$ a finite-dimensional linear space (other
variants: $\H$, $\V$), and by $\E^*$ its {\em dual} space,
composed by linear functions on $\E$. The value of function
$s \in \E^*$ at point $x \in \E$ is denoted by $\la s, x
\ra$. This notation is the same for all spaces in use.

For an operator $A: \E \to \H^*$ we denote by $A^*$ the
corresponding {\em adjoint} operator:
$$
\ba{rcl}
\la A x, y \ra & = & \la A^* y, x \ra,\quad x \in \E,\; y
\in \H.
\ea
$$
Thus, $A^*: \H \to \E^*$. A self-adjoint positive-definite
operator $B: \E \to \E^*$ (notation $B \succ 0$) defines the
Euclidean norms for the primal and dual spaces:
$$
\ba{rcl}
\| x \|_B & = & \la B x, x \ra^{1/2},\; x\in \E, \quad \| s
\|_B \; = \; \la s, B^{-1} s \ra^{1/2}, \; s \in \E^*.
\ea
$$
The sense of this notation is determined by the space of
arguments. We use the following notation for ellipsoids in
$\E$:
$$
\ba{rcl}
{\cal E}_B(x,r) & = & \{ u \in \E:\; \| u - x \|_B \leq r
\}.
\ea
$$
If in this notation parameter $r$ is missing, this means
that $r=1$.

In what follows, we often use the following simple result
from Linear Algebra. Let self-adjoint linear operators
$\Delta$ and $B$ map $\E$ to $\E^*$, and $B \succ 0$.
Then, for a tolerance parameter $\tau
> 0$ we have
\beq\label{eq-LA}
\pm \Delta \preceq \tau B \quad \Leftrightarrow \quad
\Delta B^{-1} \Delta \preceq \tau^2 B.
\eeq

For future reference, let us recall some facts from
the theory of self-concordant functions. Most of these
results can be found in Section 4 in \cite{Intro}. We use
the following notation for {\em gradient} and {\em
Hessian} of function $\Phi$:
$$
\nabla \Phi(x) \in \E^*,\quad \nabla^2 \Phi(x) \cdot h \in
\E^*,
\quad x, h \in \E.
$$

Let $\Phi$ be a self-concordant function defined on the
interior of a convex set $Q \subset \E$:
\beq\label{def-SCF}
\ba{rcl}
\left|\nabla^3 \Phi(x)[h,h,h]\right| & \leq & 2 \la \nabla^2 \Phi(x)
h,h \ra^{3/2}, \quad x \in \inter Q, \; h \in \E,
\ea
\eeq
where $\nabla^3 \Phi(x)[h_1,h_2,h_3]$ is the third
differential of function $\Phi$ at point $x$ along the
corresponding directions $h_1,h_2,h_3$. Note that $\nabla^3
\Phi(x)[h_1,h_2,h_3]$ is a trilinear symmetric form. Thus,
$$
\ba{rcl}
\nabla^3 \Phi(x)[h_1,h_2] & = & \nabla^3
\Phi(x)[h_2,h_1]\; \in \; \E^*,
\ea
$$
and $\nabla^3 \Phi(x)[h]$ is a self-adjoint linear
operator from $\E$ to $\E^*$.

Assume that $Q$ contains no straight line. Then $\nabla^2
\Phi(u)$ is nondegenerate for every $u \in \inter Q$.
Self-concordant function $\Phi$ is called
$\nu$-self-concordant barrier if
\beq\label{eq-NDec1}
\ba{rcl}
\la \nabla \Phi(u), [\nabla^2 \Phi(u)]^{-1} \nabla \Phi(u)
\ra & \leq & \nu.
\ea
\eeq
For local norms related to self-concordant functions we
use the following concise notation:
$$
\ba{rcl}
\| h \|_u & = & \la \nabla^2 \Phi(u)h,h \ra^{1/2}, \quad h
\in \E,\\
\\
\| s \|_u & = & \la s, [\nabla^2 \Phi(u)]^{-1}s \ra^{1/2},
\quad s \in \E^*.
\ea
$$
Thus, inequality (\ref{eq-NDec1}) can be written as $\|
\nabla \Phi(u) \|_u^2 \leq \nu$.

The following result is very useful.
\BT (Theorem on Recession Direction; see Section 4 in \cite{Intro} for the proof.)
If $h$ is a recession direction of the set $Q$ and $u \in
\inter Q$, then
\beq\label{eq-RDir}
\ba{rcl}
\| h \|_u & \leq & \la - \nabla \Phi(u), h \ra.
\ea
\eeq
\ET

For $u \in \inter Q$, define the Dikin ellipsoid $W_r(u)
\Def {\cal E}_{\nabla^2 \Phi(u)}(u,r)$. Then $W_r(u)
\subseteq Q$ for all $r \in [0,1)$. If $v \in W_r(u)$,
then
\beq\label{eq-GDiff}
\ba{rcl}
\la \nabla \Phi(v) - \nabla \Phi(u), v - u \ra & \geq &
{r^2 \over 1 + r}, \quad r \geq 0.
\ea
\eeq
For $r \in [0,1)$ we have
\beq\label{eq-Dik}
\ba{rcl}
(1-r)^2 \nabla^2 \Phi(u) & \preceq & \nabla^2 \Phi(v) \;
\preceq \; {1 \over (1-r)^2} \nabla^2 \Phi(u),
\ea
\eeq
\beq\label{eq-GradDiff}
\ba{rcl}
\| \nabla \Phi(v) - \nabla \Phi(u) \|_u & \leq & {r \over
1 - r},
\ea
\eeq
\beq\label{eq-HessDiff}
\ba{rcl}
\| \nabla \Phi(v) - \nabla \Phi(u) - \nabla^2 \Phi(u)(v-u)
\|_u & \leq & {r^2 \over 1 - r}.
\ea
\eeq

Finally, we need several statements on barriers for {\em
convex cones}. We call cone $K \subset \E$ {\em regular}, if
it is a closed, convex, and pointed cone with nonempty
interior. Sometimes it is convenient to write inclusion $x
\in K$ in the form $x \succeq_K 0$.

If $K$ is regular, then the {\em dual cone}
$$
\ba{rcl}
K^* & = & \left\{ s \in \E^*: \; \la s, x \ra \geq 0, \;
\forall x \in K \; \right\},
\ea
$$
is also regular. For cone $K$, we assume available a
$\nu$-normal barrier $F(x)$. This means that $F$ is
self-concordant and $\nu$-logarithmically homogeneous:
\beq\label{def-LH}
\ba{rcl}
F(\tau x ) & = & F(x) - \nu \ln \tau, \quad x \in \inter
K, \; \tau > 0.
\ea
\eeq
Note that $-\nabla F(x) \in \inter K^*$ for every $x \in \inter
K$. Equality (\ref{def-LH}) leads to many interesting
identities:
\bea
\nabla F(\tau x ) & = & \tau^{-1} \cdot \nabla F(x),
\label{eq-GT}\\
\nabla^2 F(\tau x ) & = & \tau^{-2} \cdot \nabla^2 F(x),
\label{eq-HT}\\
\la \nabla F(x), x \ra & = & - \nu, \label{eq-GX}\\
\nabla^2 F(x) \cdot x & = & - \nabla F(x), \label{eq-HX}\\
\nabla^3 F(x) [x] & = & - 2 \nabla^2 F(x), \label{eq-3X}\\
\| \nabla F(x) \|_x^2 & = & \nu, \label{eq-NDec}
\eea
where $x \in \inter K$ and $\tau > 0$. Note that the {\em
dual barrier}
$$
\ba{rcl}
F_*(s) & = & \max\limits_{x \in \inter K} \left\{\; - \la
s, x \ra - F(x) \; \right\}
\ea
$$
is a $\nu$-normal barrier for cone $K^*$. The differential
characteristics of the primal and dual barriers are
related as follows:
\beq\label{eq-PDGH}
\ba{rclcrcl}
\nabla F(- \nabla F_*(s)) & = & - s, & \quad & \nabla^2
F(-
\nabla F_*(s)) & = & [\nabla^2 F_*(s)]^{-1},\\
\\
\nabla F_*(- \nabla F(x)) & = & - x, & \quad & \nabla^2
F_*(- \nabla F(x)) & = & [\nabla^2 F(x)]^{-1},
\ea
\eeq
where $x \in \inter K$ and $s \in \inter K^*$.

For normal barriers, the Theorem on Recession Direction
(\ref{eq-RDir}) can be written both in primal and dual
forms:
\bea
\| u \|_x & \leq & \la - \nabla F(x), u \ra, \quad x \in
\inter K, u \in K,\label{eq-TRDP}\\
\| s \|_x & \leq & \la s, x \ra, \quad x \in \inter K, s
\in K^*. \label{eq-TRDD}
\eea

The following statement is very useful.
\BL\label{lm-Ell}
Let $F$ be a $\nu$-normal barrier for $K$ and $H: \E \to
\E^*$, $H \succ 0$. Assume that ${\cal E}_H (u) \subset
K$, and for some $x \in \inter K$ we have
$$
\ba{rcl}
\la \nabla F(x), u - x \ra & \geq & 0.
\ea
$$
Then, $H \succeq {1 \over 4 \nu^2} \nabla^2 F(x)$.
\EL
\proof
Let us fix an arbitrary direction $h \in \E^*$. We can assume
that
\beq\label{eq-HSign}
\ba{rcl}
\la \nabla F(x), H^{-1} h \ra & \geq & 0,
\ea
\eeq
(otherwise, multiply $h$ by $-1$). Denote $y = u + {H^{-1}
h \over \| h \|_H}$. Then $y \in K$. Therefore,
$$
\ba{rcl}
{\| H^{-1} h \|_x \over \| h \|_H} & \leq & \| u \|_x +
\|y\|_x \; \stackrel{(\ref{eq-TRDP})}{\leq} \; \la -
\nabla F(x), u \ra + \la - \nabla F(x), y \ra  \\
\\
& \stackrel{(\ref{eq-HSign})}{\leq} &  2 \la - \nabla
F(x), u \ra \; \leq \; 2 \la -\nabla F(x), x \ra \;
\stackrel{(\ref{eq-GX})}{=} \; 2 \nu.
\ea
$$
Thus, $H^{-1} \nabla^2 F(x) H^{-1} \preceq 4 \nu^2
H^{-1}$.
\qed

\BC\label{cor-Dik}
Let $x, u \in \inter K$ and $\la \nabla F(x), u - x \ra
\geq 0$. Then $\nabla^2 F(u) \succeq {1 \over 4 \nu^2}
\nabla^2 F(x)$.
\EC

\BC\label{cor-PStep}
Let $x \in \inter K$ and $u \in K$. Then $\nabla^2 F(x+u)
\preceq 4 \nu^2 \nabla^2 F(x)$.
\EC
\proof
Denote $y = x+u \in \inter K$. Then $\la \nabla F(y), x -
y \ra = \la - \nabla F(y), u \ra \geq 0$. Hence, we can
apply Corollary~\ref{cor-Dik}.
\qed

To conclude with notation, let us introduce the following
relative measure for directions in $\E$:
\beq\label{def-SH}
\ba{rcl}
\sigma_x(h) & = & \min\limits_{\rho \geq 0} \{ \rho: \;
\rho \cdot x - h \in K  \} \; \leq \; \| h \|_x, \quad x
\in \inter K, \; h \in \E.
\ea
\eeq

\section{Predicting the optimal solution}\label{sc-Pred}
\setcounter{equation}{0}

Consider the standard conic optimization problem:
\beq\label{prob-standP}
\ba{c}
\min\limits_{x \in K} \left\{\; \la c, x \ra: \; A x = b
\; \right\},
\ea
\eeq
where $c \in \E^*$, $b \in \H^*$, $A$ is a linear transformation
from $\E$ to $\H^*$, and $K \subset \E$ is a regular cone.
The problem dual to (\ref{prob-standP}) is then
\beq\label{prob-standD}
\ba{c}
\max\limits_{s \in K^*, \; y \in \H} \left\{ \la b, y \ra:
\; s + A^* y = c \; \right\}.
\ea
\eeq
Note that the feasible points of the primal and dual
problems move in the orthogonal subspaces:
\beq\label{eq-Ort}
\ba{rcl}
\la s_1 - s_2, x_1 - x_2 \ra & = & 0
\ea
\eeq
for all $x_1, x_2  \in  {\cal F}_p \Def \{ x \in K: \; A x
= b \}$, and $s_1, s_2  \in {\cal F}_d \Def\{s \in K^*: \;
s + A^* y = c \}$.

Under the {\em strict feasibility} assumption,
\beq\label{eq-Strict}
\ba{c}
\exists\; x_0 \in \inter K, \; s_0 \in \inter K^*, \; y_0
\in \H:
\quad A x_0 = b, \quad s_0 + A^* y_0 = c,
\ea
\eeq
the optimal sets of the primal and dual problems are
nonempty and bounded, and there is no duality gap
(see for instance \cite{NN}).
Moreover, a {\em primal-dual central path} $z_{\mu} \Def
(x_{\mu},s_{\mu},y_{\mu})$:
\beq\label{def-CP}
\left.
\ba{rcl}
A x_{\mu}& = & b \\
\\
c + \mu \nabla F(x_{\mu}) & = & A^* y_{\mu}\\
\\
s_{\mu} \; = \; - \mu \nabla F(x_{\mu}) &
\stackrel{(\ref{eq-PDGH}), (\ref{eq-GT})}{\Leftrightarrow}
& x_{\mu} \; = \; - \mu \nabla F_*(s_{\mu})
\ea \right\}, \quad \mu > 0,
\eeq
is well defined. Note that
\beq\label{eq-CPGap}
\ba{rcl}
\la c, x_{\mu} \ra - \la b, y_{\mu} \ra & = & \la s_{\mu},
x_{\mu} \ra \; \stackrel{(\ref{def-CP}), (\ref{eq-GX})}{=}
\; \nu \cdot \mu.
\ea
\eeq

The majority of modern strategies for solving the
primal-dual problem pair (\ref{prob-standP}),
(\ref{prob-standD}) suggest to follow this trajectory as
$\mu \to 0$. On the one hand, it is important that $\mu$
be decreased at a linear rate to attain a polynomial-time
complexity. However, on the other hand, in a small neighborhood of the
solution, it is highly desirable to switch on a
superlinear rate. Such a possibility was already
discovered for Linear Programming problems
\cite{ZTD1992,YGTZ1993,Meh1993}. There has also been significant
progress in the case of Semidefinite Programming
\cite{KSS1998,PS1998,LSZ1998,JPS1999}. In this paper, we
study more general conic problems.

For a fast local convergence of a path-following scheme,
we need to show that the predicted point
$$
\ba{rcl}
\hat z_{\mu} & = & z_{\mu} - z'_{\mu} \cdot \mu
\ea
$$
enters a small neighborhood of the solution point
$$
\ba{rcl}
z_* & = & \lim\limits_{\mu \to 0} z_{\mu}\; = \;
(x_*,s_*,y_*).
\ea
$$
It is more convenient to analyze this situation by looking at
$y$-component of the central path.

Note that $s$-component of the dual problem
(\ref{prob-standD}) can be easily eliminated:
$$
\ba{rcl}
s & = & s(y) \; \Def \; c - A^* y.
\ea
$$
Then, the remaining part of the dual problem can be
written in a more concise full-dimensional form:
\beq\label{prob-Y}
\ba{rcl}
f^* & \Def & \max\limits_{y \in \H} \{ \; \la b, y \ra: \;
y \in Q
\},\\
\\
Q & \Def & \{ y \in \H:\; c - A^* y \in K^* \}.
\ea
\eeq
In view of Assumption (\ref{eq-Strict}), interior of the
set $Q$ is nonempty. Moreover, for this set we have a
$\nu$-self-concordant barrier
$$
\ba{rcl}
f(y) & = & F_*(c - A^* y), \quad y \in \inter Q.
\ea
$$
Since the optimal set of problem (\ref{prob-Y}) is
bounded, $Q$ contains no straight line. Thus, this barrier
has a nondegenerate Hessian at every strictly feasible point.
Note that Assumption (\ref{eq-Strict}) also implies
that the linear transformation $A$ is surjective.

It is clear that $y$-component of the primal-dual central
path $z_{\mu}$ coincides with the central path of the
problem (\ref{prob-Y}):
\beq\label{eq-CPY}
\ba{rcl}
b & = & \mu \nabla f(y_{\mu}) \; = \; - \mu A \nabla F_*(c
- A^*
y_{\mu})\\
\\
& = & - \mu A \nabla F_*(s_{\mu}) \; \refEQ{def-CP} \; A
x_{\mu},
\quad \mu > 0.
\ea
\eeq

Let us estimate the quality of the following prediction
point:
$$
\ba{rcl}
p(y) & \Def & y + v(y), \quad y \in \inter Q,\\
\\
v(y) & \Def & [\nabla^2 f(y)]^{-1} \nabla f(y), \quad
s_p(y) \; \Def \; s(y) - A^* v(y).
\ea
$$
Definition of the displacement $v(y)$ is motivated by
identity (\ref{eq-HX}), which is valid for arbitrary convex
cones. Indeed, in a neighborhood of a suitably non-degenerate
solution, the barrier function should be close to the
barrier of a tangent cone centered at the solution. Hence,
the relation (\ref{eq-HX}) should be satisfied with a
reasonably high accuracy.
For every $y \in \inter Q$, we have
$$
\ba{rcl}
p(y) & = & [\nabla^2 f(y)]^{-1} \cdot \left[\nabla^2 f(y)
y + \nabla f(y)
\right]\\
\\
& = & [\nabla^2 f(y)]^{-1} \cdot \left[ A \nabla^2 F_*(c -
A^* y) A^* y
\; - \; A \nabla F_*(c - A^*y) \right]\\
\\
& \stackrel{(\ref{eq-HX})}{=} & [\nabla^2 f(y)]^{-1} \cdot
\left[ A \nabla^2 F_*(c - A^* y) A^* y
\; + \; A \nabla^2 F_*(c - A^*y)(c - A^* y) \right]\\
\\
& = & [\nabla^2 f(y)]^{-1} A \nabla^2 F_*(c - A^* y) \cdot
c.
\ea
$$
Let us choose an arbitrary pair $(s_*,y_*)$ from the
optimal solution set of the problem (\ref{prob-standD}). Then,
$$
\ba{rcl}
c & = & A^* y_* + s_*.
\ea
$$
Thus, we have proved the following representation.
\BL\label{lm-p}
For every $y \in \inter Q$ and every optimal pair
$(s_*,y_*)$ of dual problem (\ref{prob-standD}), we
have\footnote{In fact, in representation (\ref{eq-p}) we
can replace the pair $(s_*, y_*)$ by any pairs $(\bar s,
\bar y)$, satisfying condition $c = A^* \bar y + \bar s$.
However, in our analysis we are interested only in
predicting the optimal solutions.}
\beq\label{eq-p}
\ba{rcl}
p(y) & = & y_* + [\nabla^2 f(y)]^{-1} A \nabla^2 F_*(s(y))
s_*.
\ea
\eeq
\EL

\BR
Note that the right-hand side of equation (\ref{eq-p}) has
a gradient interpretation. Indeed, let us fix some $s \in
K^*$ and define the function
$$
\ba{rcl}
\phi_s(y) & = & - \la s, \nabla F_*(c - A^* y) \ra, \quad y
\in Q.
\ea
$$
Then $\nabla \phi_s(y) = A \nabla^2 F(c - A^* y) \cdot s$,
and, for self-scaled barriers $\phi_s(\cdot)$ is convex
(as well as for the barriers with negative curvature, see
Section \ref{sc-Rec}) . Thus, the representation
(\ref{eq-p}) can be rewritten as follows:
\beq\label{eq-p1}
\ba{rcl}
p(y) & = & y_* + [\nabla^2 f(y)]^{-1} \nabla
\phi_{s_*}(y).
\ea
\eeq
Note that $[\nabla^2
f(y)]^{-1}$ in the limit
acts as a projector onto
the tangent subspace to the feasible set at the solution.
\ER

To conclude this section, let us describe our prediction
abilities from the points of the central path. For that,
we need to compute derivatives of the trajectory
(\ref{def-CP}). Differentiating the last line of this
definition in $\mu$, we obtain
$$
\ba{rcl}
0 & = & - A \nabla F_*(s_{\mu}) + \mu A \nabla^2
F_*(s_{\mu}) A^* y'_{\mu}.
\ea
$$
Thus,
\beq\label{eq-DY}
\ba{rcl}
y'_{\mu} & = & {1 \over \mu} [ A \nabla^2 F_*(s_{\mu})
A^*]^{-1} A \nabla F_*(s_{\mu}) \; = \; -{1 \over \mu}
[\nabla^2f(y_{\mu})]^{-1} \nabla f(y_{\mu}).
\ea
\eeq
Therefore, we have the following representation of the
prediction point:
\beq\label{eq-PredY}
\ba{rcl}
p(y_{\mu}) & = & y_{\mu} - \mu y'_{\mu}.
\ea
\eeq
Hence, for the points of the central path, identity
(\ref{eq-p}) can be written in the following form:
\beq\label{eq-DiffY}
\ba{rcl}
y_{\mu} - \mu y'_{\mu} - y^* & = & [\nabla^2
f(y_{\mu})]^{-1} A \nabla^2 F_*(s_{\mu}) s_*.
\ea
\eeq

For the primal trajectory, we have
$$
\ba{rcl}
x'_{\mu} & \refEQ{def-CP} & {1 \over \mu} x_{\mu} + \mu
\nabla^2 F_*(s_{\mu}) A^* y'_{\mu}\; \refEQ{eq-DY} \; {1
\over \mu} x_{\mu} + \nabla^2 F_*(s_{\mu}) A^*  [ A
\nabla^2 F_*(s_{\mu}) A^*]^{-1} A \nabla F_*(s_{\mu})\\
\\
& \refEQ{eq-HX} & {1 \over \mu} x_{\mu} - \nabla^2
F_*(s_{\mu}) A^*  [ A \nabla^2 F_*(s_{\mu}) A^*]^{-1} A
\nabla^2 F_*(s_{\mu})(c - A^*y_{\mu})\\
\\
& = & {1 \over \mu} x_{\mu} - \nabla^2 F_*(s_{\mu}) A^*  [
A \nabla^2 F_*(s_{\mu}) A^*]^{-1} A \nabla^2
F_*(s_{\mu})(s_*+ A^*(y^*-y_{\mu}))\\
\\
& = & {1 \over \mu} x_{\mu} - \nabla^2 F_*(s_{\mu})
\left( s_{\mu} - s_* + A^* [ A \nabla^2 F_*(s_{\mu}) A^*]^{-1} A \nabla^2
F_*(s_{\mu})s_* \right).
\ea
$$
Using now identity (\ref{eq-HX}), definition
(\ref{def-CP}), and equation (\ref{eq-DiffY}), we obtain
the following representation:
\beq\label{eq-DX}
\ba{rcl}
x'_{\mu} & = & \nabla^2 F_*(s_{\mu})s_* - \nabla^2
F_*(s_{\mu}) A^* \left( y_{\mu} - \mu y'_{\mu} - y^*
\right).
\ea
\eeq
Due to the primal-dual symmetry of our set-up, we have the following elegant
geometric interpretation.  Let $H \Def \nabla^2 F_*(s_{\mu})$, then we have
\begin{eqnarray*}
H^{-1/2}\left(x_*-x_{\mu}+\mu x'_{\mu} \right)
& = & \mbox{ projection of $H^{-1/2}x_*$ onto the kernel of $AH^{1/2}$,}\\
H^{1/2}A^*\left(y_*-y_{\mu}+\mu y'_{\mu} \right)
& = & \mbox{ projection of $H^{1/2}s_*$ onto the image of $\left(AH^{1/2}\right)^*$.}
\end{eqnarray*}

In Section \ref{sc-MainAss}, we introduce two conditions,
ensuring the boundedness  of the derivative $x'_{\mu}$ and
a high quality of the prediction from the central dual
trajectory: $y_{\mu} - \mu y'_{\mu} - y^* \approx
O(\mu^2)$. However, first we need to decide on the way
of measuring the distances in primal and dual spaces.

\section{Measuring the distances}\label{sc-Meas}
\SetEQ

Recall that the global complexity analysis of
interior-point methods is done in an affine-invariant
framework. However, for analyzing the local convergence of
these schemes, we need to fix some Euclidean norms in the
primal and dual spaces. Recall the definitions of
Euclidean norms based on a positive definite operator $B:
\E \to \E^*$.
\beq\label{eq-B}
\ba{rcl}
\| x \|_B & = & \la B x, x \ra^{1/2}, \; x \in \E, \quad \|
s \|_B \; = \; \la s, B^{-1} s \ra^{1/2}, \; s \in \E^*.
\ea
\eeq
Using this operator, we can define operator
\beq\label{eq-G}
\ba{rcl}
G & \Def & A B^{-1} A^* : \; \H \to \H^*.
\ea
\eeq
By a Schur complement argument and the fact that $A$ is
surjective, we conclude that
\beq\label{eq-GB}
\ba{rcl}
A^* G^{-1} A & \preceq & B.
\ea
\eeq
It is convenient to choose $B$ related in a certain way to
our cones and barriers.  Sometimes it is useful to have
$B$ such that
\beq\label{eq-Bmap}
B K \subseteq K^*.
\eeq
\BL
\label{lm-B}
Let $B$ be as above, suppose $B$ satisfies (\ref{eq-Bmap})
and $u \succeq_K \pm v$.  Then $\|v\|_B \leq \|u\|_B$.
\EL
\proof
Let $B$, $u$ and $v$ be as above.  Then,
$\iprod{Bu}{u}-\iprod{Bv}{v}=\iprod{B(u-v)}{u+v}
\stackrel{(\ref{eq-Bmap})}{\geq} 0.$
\qed

\BR
As a side remark, we will see that, if $F$ has negative
curvature (see Section \ref{sc-Rec}), then $B$ chosen in
(\ref{eq-BCompNew}) satisfies (\ref{eq-Bmap}).  In the
general case, it is also possible to satisfy
(\ref{eq-Bmap}) by choosing
\[
B=\nabla^2 F(\bar x) +\nabla F(\bar x) \left[\nabla F(\bar
x)\right]^*,
\]
where $\bar x$ is an arbitrary point from $\inter K$. This
operator acts as
\[
Bh = \nabla^2 F (\bar x)h + \iprod{\nabla F(\bar x)}{h}
\cdot \nabla F(\bar x), \,\,\,\, \forall h \in \E.
\]
Then, the property (\ref{eq-Bmap}) easily follows from the
Theorem on Recession Direction. Note that
\[
\nabla^2F(\bar x) \; \preceq \; B
\stackrel{(\ref{eq-NDec})}{\preceq} (\nu+1) \cdot
\nabla^2F(\bar x).
\]
\ER

In what follows, we choose $B$ related to the primal
central path. Let us define
\beq\label{eq-BCompNew}
\ba{rcl}
B & = & \nabla^2 F(x_{\mu}) \quad \mbox{with} \quad \mu =
1.
\ea
\eeq
We need some bounds for the points of the primal and dual
central paths. Denote by $X_* \subset \E $ the set of
limit points of the primal central path, and by $S_*
\subset \E^*$ the set of limit points of the dual
central path\emph{.}
\BL\label{lm-SetBound}
If $\mu_1 \in (0,\mu_0]$, then
\beq\label{eq-BNextCP}
\ba{rcl}
\| x_{\mu_1} \|_{x_{\mu_0}} & \leq & \nu, \quad \|
s_{\mu_1} \|_{s_{\mu_0}} \; \leq \; \nu.
\ea
\eeq
In particular, for every $x_* \in X_*$ and every $s_* \in S_*$
we have:
\beq\label{eq-SetBoundNew}
\ba{rcl}
\| x_* \|_B & \leq & \nu, \quad \| s_* \|_B \; \leq \;
\nu.
\ea
\eeq
Moreover, if $\mu \in (0,1]$, then
\beq\label{eq-NHess}
\ba{rcl}
{1 \over 4 \nu^2} B & \preceq & \nabla^2 F(x_{\mu}) \;
\preceq \; {4 \nu^2 \over \mu^2} B,\\
\\
{1 \over 4 \nu^2} B^{-1} & \preceq & \nabla^2 F_*(s_{\mu})
\; \preceq \; {4 \nu^2 \over \mu^2} B^{-1}.
\ea
\eeq
\EL
\proof
Indeed,
$$
\ba{rcl}
\| x_{\mu_1} \|_{x_{\mu_0}}^2 &
\stackrel{(\ref{eq-TRDP})}{\leq} & \la - \nabla
F(x_{\mu_0}), x_{\mu_1} \ra^2 \; = \; {1 \over \mu_0^2}
\la s_{\mu_0}, x_{\mu_1} \ra^2\;
\stackrel{(\ref{def-CP})}{=} \; {1 \over \mu_0^2} [ \la c,
x_{\mu_1} \ra - \la b, y_{\mu_0} \ra ]^2 \;
\stackrel{(\ref{eq-CPGap})}{\leq} \; \nu^2.
\ea
$$
The last inequality also uses the fact that $\la c,
x_{\mu_1} \ra \leq \la c, x_{\mu_0}\ra$. Applying this
inequality with $\mu_0=1$ and taking the limit as $\mu_1
\to 0$, we obtain (\ref{eq-SetBoundNew}) in view of the
choice (\ref{eq-BCompNew}). The reasoning for the dual
central path is the same.

Further, $\la \nabla F(x_1), x_{\mu} - x_1 \ra
\stackrel{(\ref{def-CP})}{=}  \la c, x_1 - x_{\mu} \ra
\geq 0$. Therefore, applying Corollary \ref{cor-Dik}, we
get the first relation in the first line of
(\ref{eq-NHess}). Similarly, we justify the first relation
in the second line of~(\ref{eq-NHess}). Finally,
$$
\ba{rcl}
\nabla^2 F(x_{\mu}) & \refEQ{def-CP} & \nabla^2 F(- \mu
\nabla F_*(s_{\mu})) \; \refEQ{eq-HT} \; {1 \over \mu^2}
\nabla^2 F(-\nabla F_*(s_{\mu}))\; \refEQ{eq-PDGH} \; {1
\over \mu^2} [\nabla^2 F_*(s_{\mu})]^{-1}.
\ea
$$
It remains to use the first relation in the second line of
\refEQ{eq-NHess}. Using this, with the first relation, in the
second line of \refEQ{eq-NHess} and the L\"{o}wner order reversing
property of the inverse, we conclude the second relation in the
first line of \refEQ{eq-NHess}. The last unproved relation can be
justified in a similar way.
\qed

\BC
For every $\mu \in (0,1]$ we have
\beq\label{eq-D2Bound}
\ba{rcl}
\| \nabla^2 F_*(s_{\mu}) \|_B & \leq & {4 \nu^2 \over
\mu^2}.
\ea
\eeq
\EC
\proof
Indeed, for every $h \in \E^*$, we have
$$
\ba{rcl}
\| \nabla^2 F_*(s_{\mu}) h \|^2_B & = & \la B \nabla^2
F_*(s_{\mu}) h, \nabla^2 F_*(s_{\mu}) h \ra \\
\\
& \refLE{eq-NHess} & {4 \nu^2 \over \mu^2} \la h, \nabla^2
F_*(s_{\mu}) h \ra \; \refLE{eq-NHess} \; {16 \nu^4 \over
\mu^4} \la h,B^{-1} h \ra.
\ea
$$
\qed
Finally, we need to estimate the norms of the initial
data.
\BL\label{lm-A}
We have
\beq\label{eq-ANew}
\ba{rcl}
\| A \|_{G,B} &\Def & \max\limits_{h \in \E} \left\{ \| A
h
\|_G:\; \| h \|_B = 1 \right\} \; \leq \; 1,\\
\\
\| A^* \|_{B,G} &\Def & \max\limits_{y \in \H} \left\{ \|
A^* y
\|_B:\; \| y \|_G = 1 \right\} \; \leq \; 1,\\
\\
\| b \|_G & \leq & \nu^{1/2}.
\ea
\eeq
\EL
\proof
Indeed, for every $h \in \E$, we have
$$
\ba{rcl}
\| A h \|_{G,B}^2 & = & \la A h, G^{-1} A h \ra \; = \;
\max\limits_{y \in \H} \left[ 2 \la A h, y \ra - \la G y,
y
\ra \right]\\
\\
& = & \max\limits_{y \in \H} \left[ 2 \la A^*y, h \ra -
\la
A^* y, B^{-1} A^* y \ra \right] \\
\\
& \leq & \max\limits_{s \in \E^*} \left[ 2 \la s, h \ra -
\la s, B^{-1} s \ra \right] \; = \; \| h \|_B^2.
\ea
$$
Further,
$$
\ba{rcl}
\| A^* \|_{B,G} & = & \max\limits_{h \in \E, y \in \H}
\left\{ \la A^* y , h \ra:\; \| h \|_B = 1,\;
\| y \|_G = 1 \right\}\\
\\
& = & \max\limits_{h \in \E, y \in \H}
\left\{ \la Ah, y \ra:\; \| h \|_B = 1,\;
\| y \|_G = 1 \right\} \; = \; \| A \|_{G,B} \; \leq \; 1.
\ea
$$

To justify the remaining inequality, note that
$$
\ba{rcl}
\| b \|^2_G & = & \la b, G^{-1} b \ra \; = \;
\max\limits_{y \in \H} \left[ 2\la b, y \ra - \la A^* y,
B^{-1} A^* y \ra \right]\\
\\
& = & \max\limits_{y \in \H} \left[ 2\la A^* y, x_1 \ra -
\la A^* y, B^{-1} A^* y \ra \right] \\
\\
& \leq & \max\limits_{s \in \E^*} \left[ 2\la s, x_1 \ra -
\la s, B^{-1} s \ra \right] \; = \; \la B x_1, x_1 \ra
\;\\
\\
& \stackrel{(\ref{eq-BCompNew})}{= } &  \la
\nabla^2F(x_1)x_1, x_1 \ra \;
\stackrel{(\ref{eq-GX}),(\ref{eq-HX})}{=}\; \nu.
\ea
$$
\qed

\section{Main assumptions}\label{sc-MainAss}
\setcounter{equation}{0}

Now, we can introduce our main assumptions.
\BAS\label{ass-Sharp}
There exists a constant $\gamma_d > 0$ such that
\beq\label{eq-Sharp}
\ba{rcl}
f^* - \la b, y \ra \; = \; \la s, x_* \ra & \geq &
\gamma_d \| y - y_* \|_G \; \refEQ{eq-G} \; \gamma_d \| s
- s_* \|_B,
\ea
\eeq
for every $y \in Q$ (that is $s=s(y) \in {\cal F}_d$).
\EAS

Thus, we assume that the dual problem (\ref{prob-standD})
admits a {\em sharp} optimal solution. Let us derive from
Assumption \ref{ass-Sharp} that $[\nabla^2 f(y)]^{-1}$
becomes small in norm as $y$ approaches $y_*$.
\BL\label{lm-Hess}
For every $y \in \inter Q$, we have
\beq\label{eq-Hess}
\ba{rcl}
[\nabla^2 f(y)]^{-1} & \preceq & {4 \over \gamma^2_d} [f^*
- \la b, y \ra]^2 \cdot G^{-1}.
\ea
\eeq
\EL
\proof
Let us fix some $y \in \inter Q$. Consider an arbitrary
direction $h \in \H^*$. Without loss of generality, we may
assume that $\la b, [\nabla^2 f(y)]^{-1} h \ra \geq 0$
(otherwise, we can consider direction~$-h$). Since $f$ is
a self-concordant barrier, the point
$$
\ba{rcl}
y_h \; \Def \; y + {[\nabla^2 f(y)]^{-1} h \over \la h,
[\nabla^2 f(y)]^{-1} h \ra^{1/2}}
\ea
$$
belongs to the set $Q$. Therefore, in view of inequality
(\ref{eq-Sharp}), we have
$$
\ba{rcl}
\gamma_d \| y_h - y_* \|_G & \leq & f^* - \la b, y_h \ra
\; \leq \; f^* - \la b, y \ra.
\ea
$$
Hence,
$$
\ba{rcl}
{1 \over \gamma_d} [f^* - \la b, y \ra] & \geq & {\|
[\nabla^2 f(y)]^{-1} h \|_G \over \la
h, [\nabla^2 f(y)]^{-1} h \ra^{1/2}} - \| y - y_* \|_G \\
\\
& \stackrel{(\ref{eq-Sharp})}{\geq} & {\| [\nabla^2
f(y)]^{-1} h \|_G \over \la
 h, [\nabla^2 f(y)]^{-1} h \ra^{1/2}} - {1 \over \gamma_d}
[f^* - \la b, y \ra].
\ea
$$
Thus, for every $h \in \H^*$ we have
$$
\ba{rcl}
\| [\nabla^2 f(y)]^{-1} h \|_G^2 & \leq & {4 \over
\gamma^2_d} [f^* - \la b, y \ra]^2 \cdot \la h, [\nabla^2
f(y)]^{-1} h \ra.
\ea
$$
This means that
$$
\ba{rcl}
[\nabla^2 f(y)]^{-1} G [\nabla^2 f(y)]^{-1} & \preceq & {4
\over \gamma^2_d} [f^* - \la b, y \ra]^2 [\nabla^2
f(y)]^{-1},
\ea
$$
and (\ref{eq-Hess}) follows.
\qed

Now we can estimate the size of the Hessian $\nabla^2
f(y)$ with respect to the norm induced by $G$:
\beann
\| [\nabla^2 f(y)]^{-1} \|_G & \Def & \max\limits_{h \in
\H^* } \left\{ \| [\nabla^2 f(y)]^{-1} h \|_G : \; \| h
\|_G = 1 \right\}.
\eeann

\BC\label{cor-HEst}
For every $y \in \inter Q$, we have
\beq\label{eq-HEst}
\ba{rcl}
\| [\nabla^2 f(y)]^{-1} \|_G
& \leq & {4 \over \gamma^2_d} [f^* - \la b, y \ra]^2.
\ea
\eeq
Therefore, $\| v(y) \|_G \leq {2 \nu^{1/2} \over \gamma_d}
[f^* - \la b, y \ra]$.
\EC
\proof
Note that
$$
\ba{rcl}
\| [\nabla^2 f(y)]^{-1} h \|_G^2 & = & \la  h, [\nabla^2
f(y)]^{-1}G [\nabla^2 f(y)]^{-1} h \ra, \quad h \in \H^*.
\ea
$$
Hence, (\ref{eq-HEst}) follows directly from
(\ref{eq-Hess}). Further,
$$
\ba{rcl}
\| v(y) \|_G^2 & = & \la G [\nabla^2 f(y)]^{-1} \nabla
f(y), [\nabla^2 f(y)]^{-1} \nabla f(y) \ra \\
\\
& \stackrel{(\ref{eq-Hess})}{\leq} & {4 \over \gamma^2_d}
[f^* - \la b, y \ra]^2 \la  \nabla f(y), [\nabla^2
f(y)]^{-1} \nabla f(y) \ra.
\ea
$$
It remains to use inequality (\ref{eq-NDec1}).
\qed

Assumption \ref{ass-Sharp} and Lemma \ref{lm-Hess} help us
bound the norm of the right-hand side of
representation~(\ref{eq-DiffY}). However, in this
expression there is one more object, which potentially can
be large. This is the vector $\nabla^2 F_*(s_{\mu})s_*$.
Therefore, we need one more assumption.
\BAS\label{ass-Bound}
There exists a constant $\sigma_d$ such that for every
$\mu \in (0,1]$ we have
\beq\label{eq-Bound}
\ba{rcl}
\| \nabla^2 F_*(s_{\mu}) s_* \|_B \; \leq \; \sigma_d.
\ea
\eeq
\EAS

In what follows, we always suppose that Assumptions
\ref{ass-Sharp} and \ref{ass-Bound} are valid. Let us
point out their immediate consequence.
\BT\label{th-BoundCP}
For every $\mu \in (0,1]$, we have the following bounds:
\beq\label{eq-BoundYP}
\ba{rcl}
\| y_{\mu} - \mu y'_{\mu} - y_* \|_G & \leq & {4 \sigma_d
\nu^2 \over \gamma_d^2} \, \mu^2,
\ea
\eeq
\beq\label{eq-BoundXP}
\ba{rcl}
\| x'_{\mu} \|_B & \leq & \sigma_p \Def \sigma_d \left(1+
{16 \nu^4 \over \gamma_d^2} \right).
\ea
\eeq
\ET
\proof
Indeed, in view of representation (\ref{eq-DiffY}), we
have
$$
\ba{rcl}
\| y_{\mu} - \mu y'_{\mu} - y_* \|_G & \leq & \| \nabla^2
f(y_{\mu}) \|_G \cdot \| A \|_{G,B} \cdot \| \nabla^2
F_*(s_{\mu}) s_* \|_B.
\ea
$$
Thus, in view of inequalities (\ref{eq-HEst}),
(\ref{eq-ANew}), and (\ref{eq-Bound}), we have
$$
\ba{rcl}
\| y_{\mu} - \mu y'_{\mu} - y_* \|_G & \leq & {4 \sigma_d
\over \gamma_d^2} (f^* - \la b, y_{\mu} \ra)^2.
\ea
$$
Applying now identity (\ref{eq-CPGap}), we obtain
(\ref{eq-BoundYP}). Further, in view of representation
(\ref{eq-DX}), we have
$$
\ba{rcl}
\| x'_{\mu} \|_B & \leq & \| \nabla^2 F_*(s_{\mu}) s_*
\|_B + \| \nabla^2 F_*(s_{\mu}) \|_B \cdot \| A^* \|_{B,G}
\cdot \| y_{\mu} - \mu y'_{\mu} - y_* \|_G.
\ea
$$
It remains to apply inequalities (\ref{eq-Bound}),
(\ref{eq-D2Bound}), (\ref{eq-ANew}), and
(\ref{eq-BoundYP}).
\qed
\BC
There is a unique limit point of the primal central path:
$x_* = \lim\limits_{\mu \to 0} x_{\mu}$. Moreover, for every
$\mu \in (0,1]$ we have
\beq\label{eq-XRate}
\ba{rcl}
\| x_{\mu} - x_* \|_B & \leq & \sigma_p \mu.
\ea
\eeq
\EC

Note that Assumptions \ref{ass-Sharp} and \ref{ass-Bound}
do not guarantee the uniqueness of the primal optimal
solution in the problem (\ref{prob-standP}).

Let us show by some examples that our assumptions are
natural.
\BE
Consider {\em nonnegative orthant} $K = K^* = \R^n_+$ with
barriers
$$
\ba{c}
F(x) = - \sum\limits_{i=1}^n \ln x^{(i)},
\quad F_*(s) = -n - \sum\limits_{i=1}^n \ln s^{(i)}.
\ea
$$
Denote by $I_*$ the set of {\em positive} components of
the optimal dual solution $s_*$. Then denoting by $e$ the
vector of all ones, we have
$$
\ba{rcl}
\la e, \nabla^2 F_*(s_{\mu}) s_* \ra \; = \;
\sum\limits_{i \in I_*} {s_*^{(i)} \over
(s^{(i)}_{\mu})^2} & = & \sum\limits_{i \in I_*} {1 \over
s_*^{(i)}} \left( {s_*^{(i)} \over s^{(i)}_{\mu}}
\right)^2 \; \leq \; \| s_* \|_{s_{\mu}}^2 \max\limits_{i
\in I_*} {1 \over s_*^{(i)}} \;
\stackrel{(\ref{eq-BNextCP})}{\leq} \; \max\limits_{i \in
I_*} {n^2 \over s_*^{(i)}}.
\ea
$$
Since vector $\nabla^2 F_*(s_{\mu}) s_*$ is nonnegative,
we obtain for its norm an upper bound in terms of
$\max\limits_{i \in I_*} {n^2 \over s_*^{(i)}}$. Note that
this bound is valid even for degenerate dual solutions
(too many active facets in $Q$), or multiple dual optimal
solutions (which is excluded by Assumption
\ref{ass-Sharp}).

It is interesting, that we can find a bound for vector
$\nabla^2 F_*(s_{\mu}) s_*$ based on the properties of the
{\em primal} central path. Indeed, for all $i \in I_*$, we
have
$$
\ba{rcl}
\left( \nabla^2 F_*(s_{\mu}) s_*) \right)^{(i)} & = & {s_*^{(i)}
\over (s^{(i)}_{\mu})^2} \; = \; s_*^{(i)} \cdot
{(x^{(i)}_{\mu})^2 \over \mu^2} \; = \; s_*^{(i)} \cdot
\left({x^{(i)}_{\mu} - x_*^{(i)} \over \mu} \right)^2.
\ea
$$
Thus, assuming $\| x_{\mu} - x_*\|_B \leq O(\mu )$, we get
a bound for $\nabla^2 F_*(s_{\mu}) s_*$. In view of
inequality (\ref{eq-XRate}), this confirms that for Linear
Programming our assumptions are very natural.
\EE

\BE
For the cone of positive-semidefinite matrices $K = K^* =
\S^n_+$, we choose
$$
F(X) = - \ln \det X, \quad F_*(S) = -n - \ln \det S.
$$
Then,
$$
\ba{rcl}
\la I, \nabla^2 F_*(S_{\mu}) S_* \ra & = & \la I,
S_{\mu}^{-1} S_* S_{\mu}^{-1} \ra.
\ea
$$
It seems difficult to get an upper bound for this value in
terms of $\| S_* \|_{S_{\mu}}^2 = \la S_{\mu}^{-1} S_*
S_{\mu}^{-1}, S_* \ra $. However, the second approach also
works here:
$$
\ba{rcl}
\la I, S_{\mu}^{-1} S_* S_{\mu}^{-1} \ra & = & \mu^{-2}
\la X_{\mu}^2, S_* \ra \; = \; \mu^{-2} \la
(X_{\mu}-X_*)^2, S_* \ra.
\ea
$$
Thus, we get an upper bound for $ \| \nabla^2 F_*(S_{\mu})
S_* \|_B$ assuming $\| X_{\mu} - X_* \|_B \leq O(\mu)$.
\qed
\EE

Our algorithms will work with points in a small
neighborhood of the central path defined by the {\em local
gradient proximity measure}. Denote
\beq\label{eq-Prox}
\ba{rcl}
{\cal N}(\mu,\beta) & = & \left\{ y \in \H: \; \gamma(y,\mu) \Def
\| \nabla f(y) - {1 \over \mu} b \|_y \leq \beta \right\},
\quad \mu \in (0,1],\; \beta \in [0,\half].
\ea
\eeq
This proximity measure has a very familiar interpretation
in the special case of Linear Programming.  Denoting by
$S$ the diagonal matrix made up from the slack variable $s
= c-A^Ty$, notice that Dikin's affine scaling direction in
this case is given by $\left[AS^{-2}A^T\right]^{-1}b$.
Then, our predictor step corresponds to the search
direction $\left[AS^{-2}A^T\right]^{-1}A S^{-1}e$, and our
proximity measure becomes
\[
\ba{c}
\left\|AS^{-1}e - \frac{1}{\mu}
b \right\|_{AS^{-2}A^T}.
\ea
\]

Let us prove the main result of this section.
\BL\label{lm-AppBeta}
Let $y \in {\cal N}(\mu,\beta)$ with $\mu \in (0,1]$ and
$\beta \in [0,{1 \over 9}]$. Then
\beq\label{eq-AppBeta}
\ba{rcl}
\| \nabla^2 F_*(s(y)) s_* \|_B & \leq & \sigma_d + {6
\nu^2 \over \mu} \beta,
\ea
\eeq
\beq\label{eq-AppBY}
\ba{rcl}
f^* - \la b, y \ra & \leq & \kappa_1 \cdot \mu,
\ea
\eeq
where $\kappa_1 = \nu + {\beta(\beta+\sqrt{\nu}) \over 1 -
\beta}$.
\EL
\proof
Indeed,
$$
\ba{rcl}
\| s(y) - s_{\mu} \|_{s(y)} & = & \la \nabla^2 F_*(s(y))
A^*(y - y_{\mu}), A^*(y - y_{\mu}) \ra^{1/2} \; = \; \| y
- y_{\mu} \|_y \;  \Def \; r.
\ea
$$
Therefore, by (\ref{eq-Dik}) we have
$$
\ba{rcl}
(1-r)^2 \nabla^2 F_*(s_{\mu}) \; \preceq \; \nabla^2
F_*(s(y)) & \preceq & {1 \over (1 - r)^2} \nabla^2
F_*(s_{\mu}).
\ea
$$
Denote $\Delta = \nabla^2 F_*(s(y)) - \nabla^2
F_*(s_{\mu})$. Then,
$$
\ba{rcl}
\pm \Delta & \preceq & \max\left\{ {1 \over (1-r)^2} - 1,
1 - (1-r)^2 \right\} \nabla^2 F_*(s_{\mu}) \; = \; {r(2-r)
\over (1-r)^2} \nabla^2 F_*(s_{\mu}).
\ea
$$
Note that $\| \nabla^2 F_*(s(y)) s_* \|_B \refLE{eq-Bound}
\sigma_d + \| H s_* \|_B$. At the same time,
$$
\ba{rcl}
\| H s_* \|_B^2 & = & \la B H s_*, H s_* \ra \;
\refLE{eq-NHess} \; {4 \nu^2 \over \mu^2} \la [\nabla^2
F_*(s_{\mu})]^{-1} H s_*, H s_* \ra \\
\\
& \refLE{eq-LA} & {4 \nu^2 r^2 (2-r)^2 \over \mu^2
(1-r)^4} \la \nabla^2 F_*(s_{\mu})s_*, s_* \ra \;
\refLE{eq-BNextCP} \; {4 \nu^4 r^2 (2-r)^2 \over \mu^2
(1-r)^4} .
\ea
$$
Thus, $\| \nabla^2 F_*(s(y)) s_* \|_B \leq \sigma_d + {2
\nu^2 r (2-r) \over \mu (1-r)^2}$. For proving
(\ref{eq-AppBeta}), it remains to note that $r
\refLE{eq-GDiff} {\beta \over 1 -
\beta}$, and
$$
\ba{rcl}
{r(2-r) \over (1-r)^2} & = & {1 \over (1-r)^2} - 1 \; \leq
\; {(1-\beta)^2 \over (1-2\beta)^2} - 1 \; = \; {2\beta -
3 \beta^2 \over (1-2\beta)^2} \; < \; {\beta(2-3\beta)
\over 1 - 4 \beta}\; \leq \; 3
\beta, \quad \beta \in [0,{1 \over 9}].
\ea
$$

To establish (\ref{eq-AppBY}), note that
$$
\ba{rcl}
\frac{1}{\mu} \left[f^* - \la b, y \ra\right] & = & \frac{1}{\mu}
\left[\la b, y^*-y_{\mu} \ra + \la b, y_{\mu}-y \ra \right]\\
\\
& \stackrel{(\ref{eq-CPGap})}{\leq} &
\nu + \la \frac{b}{\mu}, y_{\mu}-y \ra \\
\\
& = & \nu + \la -\nabla f(y) + \frac{b}{\mu}, y_{\mu}-y \ra
+ \la \nabla f(y), y_{\mu}-y \ra \\
\\
& \leq & \nu + \| \nabla f(y) -\frac{b}{\mu}\|_y \cdot
\|y - y_{\mu}\|_y + \|\nabla f(y)\|_y \cdot \|y -y_{\mu}\|_y \\
\\
& \leq & \nu + \beta \frac{\beta}{1-\beta} +\sqrt{\nu}
\frac{\beta}{1-\beta},
\ea
$$
where the last inequality follows from the assumptions of
the lemma and (\ref{eq-NDec1}).
\qed

Now, we can put all our observations together.
\BT\label{th-DEst}
Let dual problem (\ref{prob-standD}) satisfy Assumptions
\ref{ass-Sharp} and \ref{ass-Bound}. If for some $\mu \in
(0,1]$ and $\beta \in [0,{1 \over 9}]$ we have $y \in
{\cal N}(\mu,
\beta)$, then
\beq\label{eq-DEst}
\ba{rcl}
\| p(y) - y_* \|_G & \leq & {4 \over \gamma_d^2} \left(
\sigma_d + {6 \nu^2 \over \mu} \beta \right) \la b, y -
y_* \ra^2 \; \leq \; {4 \nu \over \gamma_d^2} \left(
\sigma_d + {6 \nu^2 \over \mu} \beta \right)\cdot \| y -
y_* \|^2_G.
\ea
\eeq
\ET
\proof
Indeed, in view of representation (\ref{eq-p}), we have
$$
\ba{rcl}
\| p(y) - y_* \|_G & \leq & \| [\nabla^2 f(y)]^{-1} \|_G
\cdot \| A \|_{G,B} \cdot \| \nabla^2 F_*(s(y)) s_* \|_B.
\ea
$$
Now, we can use inequalities (\ref{eq-HEst}),
(\ref{eq-ANew}), and (\ref{eq-AppBeta}). For justifying
the second inequality, we apply the third bound in
(\ref{eq-ANew}).
\qed

The most important consequence of the estimate
(\ref{eq-DEst}) consists in the necessity to keep the
neighborhood size parameter $\beta$ at the same order
as the (central path) penalty
parameter $\mu$. In our reasoning below, we often use
\beq\label{eq-BMu}
\ba{rcl}
\beta & = & {1 \over 9} \mu.
\ea
\eeq

\section{Efficiency of the predictor step}\label{sc-PStep}
\setcounter{equation}{0}

Let us estimate now the efficiency of the predictor step
with respect to the local norm.
\BL\label{lm-Loc}
If $y \in {\cal N}(\mu,\beta)$, then
\beq\label{eq-Loc}
\ba{rcl}
\| p(y) - y_* \|_y & \leq & \kappa_2 \cdot [f^* - \la b, y
\ra ] \; \leq \; \mu \cdot \kappa,
\ea
\eeq
where $\kappa_2 = {2 \over \gamma_d} \left(\sigma_d + {6
\nu^2 \over \mu} \beta \right)$, and $\kappa = \kappa_1
\cdot \kappa_2$.
\EL
\proof
Indeed,
$$
\ba{rcl}
\| p(y) - y_* \|_y^2 & \stackrel{(\ref{eq-p})}{=} & \la A
\nabla^2 F_*(s(y)) s_*, [\nabla^2 f(y)]^{-1} A \nabla^2
F_*(s(y))
s_* \ra\\
\\
& \stackrel{(\ref{eq-Hess})}{\leq } & {4 \over \gamma_d^2}
[f^* - \la b, y \ra ]^2 \la A \nabla^2 F_*(s(y)) s_*,
G^{-1} A \nabla^2 F_*(s(y)) s_* \ra \\
\\
& \stackrel{(\ref{eq-GB})}{\leq } & {4 \over \gamma_d^2}
[f^* - \la b, y \ra ]^2 \la B \nabla^2 F_*(s(y)) s_*,
\nabla^2 F_*(s(y)) s_* \ra.
\ea
$$
It remains to use the bounds (\ref{eq-AppBeta}) and
(\ref{eq-AppBY}).
\qed

Since $\| y_* - y \|_y \geq 1$, inequality (\ref{eq-Loc})
demonstrates a significant drop in the distance to the
optimal point after a full predictor step.
The following fact is also useful.
\BL\label{lm-Null}
For every $y \in Q$, we have $A \cdot \nabla^2 F_*(s(y))
\cdot s_p(y)= 0$.
\EL
\proof
Indeed, $A \nabla^2 F_*(s(y)) s_p(y) = A \nabla^2
F_*(s(y))(s(y) - A^*v(y)) \stackrel{(\ref{eq-HX})}{=} - A
\nabla F_*(s(y)) - \nabla f(y) = 0$.
\qed

\BC
If ${\cal F}_p$ is bounded, then the point $\nabla^2
F_*(s(y)) \cdot s_p(y) \notin K$ (therefore, it
is infeasible for the primal
problem (\ref{prob-standP})).
\EC

We can show now that a large predictor step can still keep
dual feasibility. Denote
$$
\ba{rcl}
y(\alpha) & = & y + \alpha v(y), \quad \alpha \in [0,1].
\ea
$$
\BT\label{th-Feas}
Let $y \in {\cal N}(\mu,\beta)$ with $\mu \in (0,1]$ and
$\beta \in [0,{1 \over 9}]$. Then, for every $r \in
(0,1)$, the point $y(\hat \alpha)$ with
\beq\label{eq-Alpha}
\ba{rcl}
\hat \alpha & \Def & {r \over r + \kappa_2 [ f^* - \la b,
y \ra]}
\ea
\eeq
belongs to $Q$. Moreover,
\beq\label{eq-Quad}
\ba{rcl}
f^* - \la b, y(\hat \alpha) \ra & \leq & \kappa_3 \cdot
[f^* - \la b, y \ra]^2,
\ea
\eeq
where $\kappa_3 = \kappa_2 \cdot \left({1 \over r} + {2
\sqrt{\nu} \over \gamma_d} \right)$.
\ET
\proof
Consider the Dikin ellipsoid $W_r(y) = \{ u \in \H: \| u - y \|_y
\leq r\}$. Since $W_r(y) \subseteq Q$, its convex
combination with point $y_*$, defined as
$$
\ba{rcl}
Q(y) & = & \{ u \in \H: \; \| u - (1-t)y - ty_* \|_y \leq r(1-t),
\; t \in [0,1] \},
\ea
$$
is contained in $Q$. Note that
$$
\ba{rcl}
\| y(\hat \alpha) - (1-\hat \alpha)y - \hat \alpha y_*
\|_y & = & \hat \alpha \| p(y) - y_* \|_y \\
\\
& \stackrel{(\ref{eq-Loc})}{\leq} & \kappa_2 \hat \alpha
[f^* - \la b, y \ra] \; \stackrel{(\ref{eq-Alpha})}{=}\;
r(1 - \hat \alpha).
\ea
$$
Hence, $y(\hat \alpha) \in Q$. Further,
$$
\ba{rcl}
f^* - \la b, y(\hat \alpha) \ra & = & (1 - \hat \alpha)
[f^* - \la b, y \ra] + \hat \alpha \la b, y_* - p(y) \ra \\
\\
& \leq & {\kappa_2 \over r} [f^* - \la b, y \ra]^2 + \| b
\|_G \cdot \| p(y) - y_* \|_G.
\ea
$$
Since $\| b \|_G \stackrel{(\ref{eq-ANew})}{\leq}
\sqrt{\nu}$ and
$$
\ba{rcl}
\| p(y) - y_* \|_G & \stackrel{(\ref{eq-DEst})}{\leq} & {2
\kappa_2 \over \gamma_d} [ f^* - \la b, y \ra ]^2,
\ea
$$
we obtain the desired inequality (\ref{eq-Quad}).
\qed

Denote by $\bar \alpha(y)$, the maximal feasible step
along direction $v(y)$:
$$
\ba{rcl}
\bar \alpha(y) & = & \max\limits_{\alpha \geq 0} \{
\alpha: \; y + \alpha v(y) \in Q \}.
\ea
$$
Let us show that $\bar \alpha=\bar \alpha(y)$ is large
enough. In general,
\beq\label{eq-DikLow}
\ba{rcl}
\bar \alpha(y) & \geq & {1 \over \| v(y) \|_y} \;
\stackrel{(\ref{eq-NDec1})}{\geq} \; {1 \over \nu^{1/2}}.
\ea
\eeq
However, in a small neighborhood of the solution, we can
establish a better bound.
\BT\label{lm-BarAlpha} Let $y
\in {\cal N}(\mu,\beta)$ with $\mu \in (0,1]$ and $\beta
\in [0,{1 \over 9}]$. Then
\beq\label{eq-BAlpha1}
\ba{rcl}
1 - \bar \alpha(y) & \leq & { \kappa \mu \over 1 + \kappa
\mu}.
\ea
\eeq
Moreover, if $\mu$ is small enough:
\beq\label{eq-MuBound}
\ba{rcl}
\mu & < {1 - 2 \beta \over \kappa},
\ea
\eeq
then
\beq\label{eq-BAlpha11}
\ba{rcl}
\bar \alpha(y)-1 & \leq & { \kappa \mu \over 1 - \kappa
\mu -2\beta},
\ea
\eeq
and
\beq\label{eq-BAlpha2}
\ba{rcl}
\| y(\bar \alpha) - y^* \|_y & \leq & \kappa \mu \left(1 +
{\sqrt{\nu} \over 1 - \kappa \mu - 2 \beta} \right).
\ea
\eeq
\ET
\proof
Since for every $r \in (0,1)$
$$
\ba{rcl}
\bar \alpha & \stackrel{(\ref{eq-Alpha})}{\geq} & \hat \alpha = {r \over r +
\kappa_2 [f^* - \la b, y \ra]},
\ea
$$
we have $1 - \bar \alpha \leq { \kappa_2 [f^* - \la b, y
\ra] \over 1 +  \kappa_2 [f^* - \la b, y \ra]}
\stackrel{(\ref{eq-AppBY})}{\leq} { \kappa \mu \over 1 +
\kappa \mu}$, which is (\ref{eq-BAlpha1}). On the other
hand,
$$
\ba{rcl}
\la b, v(y) \ra & = & \la b, [f''(y)]^{-1} (f'(y) - {1
\over \mu} b + {1 \over \mu} b \ra \; \geq \; {1 \over
\mu} \| b \|_y^2 - \beta \| b \|_y\\
\\
& = & \| b \|_y ( \| f'(y) + {1 \over \mu}b - f'(y) \|_y -
\beta) \; \geq \; \| b \|_y ( \| f'(y) \|_y - 2 \beta)\\
\\
& = & \| b \|_y ( \| p(y) - y \|_y - 2 \beta) \; \geq \;
\| b \|_y ( \| y - y_* \|_y - \| p(y) - y_* \|_y - 2
\beta).
\ea
$$
Thus, using the estimate (\ref{eq-Loc}) and the bound $\|
y - y_* \|_y \geq 1$ (since $y_*$ is on the boundary of $Q$), we get
\beq\label{eq-BPos}
\ba{rcl}
\la b, v(y) \ra & \geq & \| b \|_y (1 - \kappa \mu - 2
\beta).
\ea
\eeq
Therefore, condition (\ref{eq-MuBound}) guarantees that
$\la b, v(y) \ra > 0$. Define $\tilde \alpha = {f^* - \la
b, y \ra \over \la b, v(y) \ra}$. Then
$$
\ba{rcl}
\la b, y + \tilde \alpha v(y) \ra & = & f^*.
\ea
$$
Therefore, $y + \tilde \alpha v(y) \not\in Q$. Hence,
$$
\ba{rcl}
\bar \alpha & \leq & \tilde \alpha \; = \; 1 + {\la b, y_*
- y - v(y) \ra \over \la b, v(y) \ra} \; \refLE{eq-Loc} \;
1 + { \| b \|_y \cdot \kappa \mu \over \la b, v(y) \ra} \;
\refLE{eq-BPos} \; 1 + {\kappa \mu \over 1 - \kappa \mu -
2 \beta}.
\ea
$$

Further,
$$
\ba{rcl}
\| y(\bar \alpha) - y^* \|_y & \leq & \| y(\bar \alpha) -
p(y) \|_y + \| p(y) - y^* \|_y \;
\stackrel{(\ref{eq-Loc})}{\leq} \; |1-\bar \alpha| \cdot
\|
v(y) \|_y + \kappa \cdot \mu\\
\\
&\stackrel{(\ref{eq-NDec1})}{\leq} & |1-\bar \alpha|
\sqrt{\nu} + \kappa \cdot \mu .
\ea
$$
Taking into account that in view of (\ref{eq-BAlpha1}) and
(\ref{eq-BAlpha11}) $|1 - \bar \alpha| \leq {\kappa \mu
\over 1 - \kappa \mu - 2 \beta}$, we get inequality
(\ref{eq-BAlpha2}).
\qed

Despite the extremely good progress in function value, we
have to worry about the distance to the central path
and we cannot yet fully appreciate the new point $y(\hat \alpha)$.
Indeed, for getting close again to the central path, we
need to find an approximate solution to the auxiliary
problem
$$
\min\limits_y \{ f(y): \la b, y \ra = \la b, y(\hat
\alpha) \ra \}.
$$
In order to estimate the complexity of this {\em corrector
stage}, we need to develop some bounds on the growth of
the gradient proximity measure.

\section{Barriers with negative curvature}\label{sc-Rec}
\setcounter{equation}{0}

\BD\label{def-NC}
Let $F$ be a normal barrier for the regular cone $K$. We say that $F$
has {\em negative curvature} if for every $x \in \inter K$
and $h \in K$ we have
\beq\label{eq-NC}
\ba{rcl}
\nabla^3F(x)[h] & \preceq & 0.
\ea
\eeq
\ED
It is clear that self-scaled barriers have negative
curvature (see \cite{NT1}). Some other important barriers,
like the negation of the logarithms of {\em hyperbolic polynomials} (see
\cite{Guler1997}) also share this property.

\BT\label{th-NC}
Let $K$ be a regular cone and $F$ be a normal barrier for
$K$. Then, the following statements are equivalent:
\begin{enumerate}
\item
$F$ has negative curvature;
\item
for every $x \in \inter K$ and $h \in \E$ we have
\beq\label{eq-NC3}
\ba{rcl}
-\nabla^3F(x)[h,h] & \in & K^*;
\ea
\eeq
\item
for every $x \in \inter K$ and for every $h \in \E$
such that $x+h \in \inter K$, we have
\beq\label{eq-NC2}
\ba{rcl}
\nabla F(x+h) - \nabla F(x) \; \preceq_{K^*} \; \nabla^2
F(x) h.
\ea
\eeq
\end{enumerate}
\ET
\proof
Let $F$ have negative curvature. Then, for every $h \in \E$
and $u \in K$ we have
\beq\label{eq-T30}
\ba{rcl}
0 & \geq & \nabla^3 F(x)[h,h,u] \; = \; \la \nabla^3
F(x)[h,h], u \ra.
\ea
\eeq
Clearly, this condition is equivalent to (\ref{eq-NC3}).
On the other hand, from (\ref{eq-NC3}) we have
$$
\ba{rcl}
\nabla F(x+h) - \nabla F(x) - \nabla^2 F(x) h & = &
\int\limits_0^1 \nabla^3 F(x + \tau h) [h, h]\, d \tau \;
\preceq_{K^*} \; 0.
\ea
$$
Note that we can replace in (\ref{eq-NC2}) $h$ by $\tau
h$, divide everything by $\tau^2$, and take the limit as
$\tau \to 0^+$. Then we arrive back at the inclusion
(\ref{eq-NC3}).
\qed

\BT\label{th-Mon}
Let the curvature of $F$ be negative. Then for every
$x \in K$, we have
\beq\label{eq-NCPos}
\ba{rcl}
\nabla^2 F(x) h & \succeq_{K^*} & 0, \quad \forall h \in
K,
\ea
\eeq
and, consequently,
\beq\label{eq-NC1}
\ba{rcl}
\nabla F(x+h) - \nabla F(x) & \succeq_{K^*} & 0.
\ea
\eeq
\ET
\proof
Let us prove that $\nabla^2 F(x) h \in K^*$ for $h \in K$.
Assume first that $h \in \inter K$. Consider the following
vector function:
$$
s(t) = \nabla^2 F(x+th) h \in \E^*, \quad t \geq 0.
$$
Note that $s'(t) = \nabla^3F(x+th)[h,h]
\stackrel{(\ref{eq-NC3})}{\preceq_{K^*}} 0$. This means
that
$$
\ba{rcl}
\nabla^2 F(x)h & \succeq_{K^*} & \nabla^2 F(x+th) h \;
\stackrel{(\ref{eq-HT})}{=} {1 \over t^2} \nabla^2 F(h+{1
\over t}x) h.
\ea
$$
Taking the limit as $t \to \infty$, we get $\nabla^2 F(x)h
\in K^*$. By continuity arguments, we can extend this
inclusion onto all $h \in K$. Therefore,
$$
\ba{rcl}
\nabla F(x+h) & = & \nabla F(x) + \int\limits_0^1 \nabla^2
F(x + \tau h) h \, d \tau \; \succeq_{K^*} \; \nabla F(x).
\ea
$$
\qed

As we have proved, if $F$ has negative curvature, then
$\nabla^2 F(x) K \subseteq K^*$, for every $x \in \inter
K$. This property implies that the situations when both
$F$ and $F_*$ has negative curvature are very seldom.
\BL\label{lm-NCS}
Let both $F$ and $F^*$ have negative curvature. Then $K$
is a symmetric cone.
\EL
\proof
Indeed, for every $x \in \inter K$ we have $\nabla^2 F(x)
K \subseteq K^*$. Denote $s = - \nabla F(x)$. Since $F_*$
has negative curvature, then $\nabla^2 F_*(s) K^*
\subseteq K$. However, since $\nabla^2 F_*(s)
\stackrel{(\ref{eq-PDGH})}{=} [\nabla^2 F(x)]^{-1}$, this
means $K^* \subseteq \nabla^2 F(x) K$. Thus $K^* =
\nabla^2 F(x) K$. Now, using the same arguments as
in~\cite{NT1} it is easy to prove that for every pair $x
\in \inter K$ and $s \in \inter K^*$ there exists a
scaling point $w \in \inter K$ such that $s = \nabla^2
F(w) x$ (this $w$ can be taken as the minimizer of the
convex function $-\la \nabla F(w), x \ra + \la s, w \ra$).
Thus, we have proved that $K$ is homogeneous and
self-dual. Hence, it is symmetric.
\qed

\BT\label{th-NCHess}
Let $K$ be a regular cone and $F$ be a normal barrier for
$K$, which has negative curvature. Further let $x, x + h
\in \inter K$. Then for every $\alpha \in [0,1)$ we have
\beq\label{eq-NCHUp}
\ba{rcl}
{1 \over (1 + \alpha \sigma_x(h))^2} \nabla^2 F(x) \;
\preceq \; \nabla^2 F(x + \alpha h) & \preceq & {1 \over
(1 - \alpha)^2} \nabla^2 F(x).
\ea
\eeq
\ET
\proof
Indeed,
$$
\ba{rcl}
\nabla^2 F(x + \alpha h) & = & \nabla^2 F((1-\alpha)x +
\alpha (x+h)) \\
\\
& \stackrel{(\ref{eq-NC})}{\preceq} & \nabla^2
F((1-\alpha)x) \; \stackrel{(\ref{eq-HT})}{=} \; {1 \over
(1-\alpha)^2} \nabla^2 F(x).
\ea
$$

Further, denote $\bar x = x - {h \over \sigma_x(h)}$. By
definition, $\bar x \in K$. Note that
$$
\ba{rcl}
x & = & (x + \alpha h) + {\alpha \sigma_x(h) \over 1 +
\alpha \sigma_x(h)} (\bar x - (x+\alpha h)).
\ea
$$
Therefore, by the second inequality in (\ref{eq-NCHUp}),
we have
$$
\ba{rcl}
\nabla^2 F(x) & \preceq & (1 + \alpha \sigma_x(h))^2
\nabla^2 F(x + \alpha h).
\ea
$$
\qed

\section{Bounding the growth of the proximity
measure}\label{sc-Growth}
\setcounter{equation}{0}

Let us analyze now our predictor step
$$
y(\alpha) = y + \alpha v(y), \quad \alpha \in [0,\bar
\alpha],
$$
where $\bar \alpha = \bar \alpha(y)$. Denote $\bar s =
s(y(\bar \alpha)) \in K^*$.
\BL
Let $F_*$ have negative curvature.  Then,
for every $\alpha \in [0, \bar{\alpha})$, we have
\beq\label{eq-BDeltaN}
\ba{rcl}
\delta_y(\alpha) & \Def & \| \nabla f(y(\alpha)) - {\bar
\alpha \over \bar \alpha - \alpha} \nabla f(y) \|_{G} \;
\leq \; {\alpha \bar \alpha \over (\bar \alpha -
\alpha)^2} \| \nabla^2 F_*(s(y)) \bar s \|_{B}.
\ea
\eeq
\EL
\proof
Indeed,
$$
\ba{rl}
& \delta_y^2(\alpha) = \la G^{-1} (\nabla f(y(\alpha)) -
{\bar \alpha \over \bar \alpha - \alpha} \nabla f(y)),
\nabla f(y(\alpha)) - {\bar \alpha \over
\bar \alpha - \alpha} \nabla f(y) \ra\\
\\
= & \la G^{-1} A(\nabla F_*(s(y(\alpha))) - {\bar \alpha
\over \bar \alpha - \alpha} \nabla F_*(s(y))), A(\nabla
F_*(s(y(\alpha))) - {\bar \alpha \over
\bar \alpha - \alpha} \nabla F_*(s(y)) ) \ra\\
\\
\stackrel{(\ref{eq-GB})}{\leq} & \la B (\nabla
F_*(s(y(\alpha))) - \nabla F_*((1 - {\alpha \over
\bar \alpha})s(y))), \nabla F_*(s(y(\alpha))) - \nabla F_*( (1 - {\alpha \over
\bar \alpha})s(y))  \ra.
\ea
$$
Note that $y(\alpha) = y + {\alpha \over \bar
\alpha}(y(\bar \alpha) - y)$. Therefore,
$$
\ba{rcl}
s(y(\alpha)) & = & (1 - {\alpha \over \bar \alpha}) s(y) +
{\alpha \over \bar \alpha}\bar s .
\ea
$$
Since $F_*$ has negative curvature, we have
\[
x' \Def \nabla F_*(s(y(\alpha))) - \nabla F_*
\left(\left(1 - {\alpha \over
\bar \alpha}\right)s(y)\right) \succeq_K 0.\]
Using (\ref{eq-NC2}) in Theorem \ref{th-NC}, we also have
\[
x'' \Def \nabla^2 F_*\left(\left(1 - {\alpha \over\bar \alpha}
\right) s(y) \right) \cdot \left({\alpha \over\bar \alpha} \bar{s}
\right) \succeq_K x'.
\]
Thus, $x'' \succeq_K \pm x'$.  By (\ref{eq-NCPos}) in Theorem \ref{th-Mon}
and
Lemma \ref{lm-B}, we obtain
$\iprod{Bx''}{x''} \geq \iprod{Bx'}{x'}$ which
gives the desired conclusion.
\qed

Note that at the predictor stage, we need to choose the
rate of decrease of the penalty parameter (central path
parameter) $\mu$ as a function
of the predictor step size $\alpha$. Inequality
(\ref{eq-BDeltaN}) suggests the following dependence:
\beq\label{eq-Mu1}
\ba{rcl}
\mu(\alpha) & \approx & \left(1 - {\alpha \over \bar
\alpha} \right) \cdot \mu.
\ea
\eeq
However, if $\bar \alpha$ is close to its lower limit
(\ref{eq-DikLow}), this strategy may be too aggressive.
Indeed, in a small neighborhood of the point $y$ we can
guarantee only
\beq\label{eq-EstLoc}
\ba{rcl}
\| \nabla f(y(\alpha))  - (1+\alpha) \nabla f(y) \|_y & =
& \| \nabla f(y(\alpha))  -  \nabla f(y) - \alpha \nabla^2
f(y) v(y) \|_y\\
\\
& \stackrel{(\ref{eq-HessDiff})}{\leq} & {\alpha^2 \| v(y)
\|_y^2 \over 1 - \alpha \| v(y) \|_y}.
\ea
\eeq
In this situation, a more reasonable strategy for decreasing
$\mu$ seems to be:
\beq\label{eq-Mu2}
\ba{rcl}
\mu(\alpha) & \approx & {\mu \over 1 + \alpha}.
\ea
\eeq
It appears that it is possible to combine both
strategies (\ref{eq-Mu1}) and (\ref{eq-Mu2})
in a single expression. Denote
$$
\ba{rcl}
\xi_{\bar \alpha}(\alpha) & = & 1 + {\alpha \bar \alpha
\over
\bar \alpha - \alpha}, \quad \alpha \in [0,\bar \alpha).
\ea
$$
Note that
\beq\label{eq-Xi}
\ba{rcl}
\xi_{\bar \alpha}(\alpha) & = & 1 + \alpha + {\alpha^2
\over
\bar \alpha - \alpha} \; = \; {\bar \alpha \over \bar
\alpha - \alpha} - {\alpha (1 - \bar \alpha) \over \bar
\alpha - \alpha}.
\ea
\eeq
Let us prove an upper bound for the growth of the local
gradient proximity measure along direction $v(y)$, when
the penalty parameter is divided by the factor $\xi_{\bar
\alpha}(\alpha)$.
\BT\label{th-NCG}
Suppose $F_*$ has negative curvature. Let $y \in {\cal
N}(\mu, \beta)$ with $\mu \in (0,1]$ and $\beta \in [0,{1
\over 9}]$, satisfying condition (\ref{eq-MuBound}). Then,
for $y(\alpha) = y + \alpha v(y)$ with $\alpha \in (0,
\bar \alpha)$, we have
\beq\label{eq-Main}
\ba{c}
\gamma\left(y(\alpha),{\mu \over \xi_{\bar
\alpha}(\alpha)}\right) \; \leq \; \Gamma_{\mu}(y,\alpha)
\Def \;
\left(1 + \alpha \cdot \sigma_{s(y)}\left(-A^*v(y)\right)
\right) \| \nabla
f(y(\alpha))  - {\xi_{\bar \alpha}(\alpha) \over \mu} \cdot b \|_y\\
\\
\; \leq \; (1 + \alpha \cdot \sigma_{s(y)}(-A^* v(y)))
\left[ \gamma_1(\alpha) +
\beta \cdot \left(1 + {\alpha \bar \alpha \over \bar
\alpha
- \alpha}\right) \right],\\
\\
\gamma_1(\alpha) \; \Def \; \| \nabla f(y(\alpha)) -
\xi_{\bar \alpha}(\alpha)
 \nabla f(y) \|_y \\
\\
\leq \; { \alpha \bar \alpha \mu \over (\bar \alpha -
\alpha)^2}
\left(  \left(1 - {\alpha \over \bar \alpha} \right)
\kappa \sqrt{\nu} + {2 \kappa_1 \over \gamma_d}
\left[ \sigma_d + {6 \nu^2 \beta \over \mu} + 2 \kappa
\nu \left(1 + {\sqrt{\nu} \over 1 - \kappa \mu - 2
\beta} \right)
{1 -
\beta \over 1 - 2
\beta} \right] \right).
\ea
\eeq
\ET
\proof
Indeed,
$$
\ba{rcl}
\gamma\left(y(\alpha),{\mu \over \xi_{\bar
\alpha}(\alpha)}\right) & = & \| \nabla
f(y(\alpha)) - {\xi_{\bar \alpha}(\alpha) \over \mu} \cdot b \|_{y(\alpha)}\\
\\
& \stackrel{(\ref{eq-NCHUp})}{\leq} &
(1 + \alpha \sigma_{s(y)}(-A^* v(y))) \cdot \| \nabla
f(y(\alpha)) -
 {\xi_{\bar \alpha}(\alpha) \over \mu} \cdot b \|_{y}.
\ea
$$
Further,
$$
\ba{rcl}
\| \nabla f(y(\alpha)) - {\xi_{\bar \alpha}(\alpha) \over
\mu} \cdot b \|_{y} & \leq & \gamma_1(\alpha) + \xi_{\bar
\alpha}(\alpha) \| \nabla f(y) - {1 \over \mu} b \|_{y}.
\ea
$$
Since $y \in {\cal N}(\mu,\beta)$, the last term does not
exceed $\beta \cdot \xi_{\bar \alpha}(\alpha)$. Let us
estimate now $\gamma_1(\alpha)$.
$$
\ba{rcl}
\gamma_1(\alpha) & \stackrel{(\ref{eq-Xi})}{\leq} & \|
\nabla f(y(\alpha)) - {\bar \alpha \over \bar \alpha -
\alpha} \nabla f(y) \|_y + {\alpha (1 - \bar \alpha) \over
\bar \alpha - \alpha} \| \nabla f(y) \|_y\\
\\
& \stackrel{(\ref{eq-BAlpha1})}{\leq} & \| \nabla
f(y(\alpha)) - {\bar \alpha \over \bar \alpha - \alpha}
\nabla f(y) \|_y + {\alpha \over \bar \alpha - \alpha}
\cdot {\kappa \mu \sqrt{\nu} \over 1 + \kappa \mu}.
\ea
$$
For the second inequality above, we also used (\ref{eq-NDec1}).
Note that
$$
\ba{rcl}
&  & \| \nabla f(y(\alpha)) - {\bar \alpha \over \bar
\alpha - \alpha} \nabla f(y)
\|_y^2\\
\\
& = & \la [\nabla^2 f(y)]^{-1}(\nabla f(y(\alpha)) - {\bar
\alpha \over \bar \alpha - \alpha} \nabla f(y)), \nabla
f(y(\alpha)) - {\bar \alpha \over \bar \alpha - \alpha}
\nabla f(y) \ra\\
\\
& \stackrel{(\ref{eq-Hess})}{\leq} & {4 \over \gamma_d^2}
[f^* - \la b, y \ra]^2 \cdot  \| \nabla f(y(\alpha)) -
{\bar \alpha \over
\bar \alpha - \alpha} \nabla f(y) \|_{G}^2 \;
\stackrel{(\ref{eq-AppBY})}{\leq} \; {4 \kappa_1^2 \mu^2
\over \gamma_d^2} \cdot \delta_y^2(\alpha).
\ea
$$
Moreover,
$$
\ba{rcl}
\delta_y(\alpha) & \stackrel{(\ref{eq-BDeltaN})}{\leq} &
{\alpha \bar \alpha \over (\bar \alpha - \alpha)^2} \|
\nabla^2 F_*(s(y))
\bar s \|_{B}\\
\\
& \leq & {\alpha \bar \alpha \over (\bar \alpha -
\alpha)^2}
\left[ \| \nabla^2 F_*(s(y))
s_* \|_{B} + \| \nabla^2 F_*(s(y)) (\bar s - s_*)
\|_{B} \right]\\
\\
& \stackrel{(\ref{eq-AppBeta})}{\leq} & {\alpha \bar
\alpha \over (\bar \alpha - \alpha)^2}
\left[ \sigma_d + {6
\nu^2 \over \mu} \beta
 + \| \nabla^2 F_*(s(y)) (\bar s - s_*)
\|_{B} \right].
\ea
$$
It remains to estimate the last term.

Denote $r = \| s(y) - s_{\mu} \|_{s(y)}
\stackrel{(\ref{eq-GDiff})}{\leq} {\beta \over 1-\beta}$.
Then
$$
\ba{rcl}
B & \refLE{eq-NHess} & {4 \nu^2 \over \mu^2} [\nabla^2
F_*(s_{\mu})]^{-1} \; \stackrel{(\ref{eq-Dik})}{\preceq}
\; {4 \nu^2 \over \mu^2 (1-r)^2} [\nabla^2
F_*(s(y))]^{-1}.
\ea
$$
Therefore,
$$
\ba{rcl}
 \| \nabla^2 F_*(s(y)) (\bar s - s_*)
\|_{B} & \leq & {2 \nu \over \mu(1-r)} \| y(\bar \alpha) -
y^* \|_y \; \stackrel{(\ref{eq-BAlpha2})}{\leq} \; 2
\kappa \nu \left(1 + {\sqrt{\nu} \over 1 - \kappa \mu - 2
\beta} \right)\cdot {1-\beta \over 1-2\beta}.
\ea
$$
Putting all the estimates together, we obtain the claimed
upper bound on $\gamma_1(\alpha)$.

\qed

Taking into account the definition of $\xi_{\bar
\alpha}(\alpha)$, we can see that our predictor-corrector
scheme with neighborhood size parameter $\beta = O(\mu)$ has local
superlinear convergence.

\section{Polynomial-time path-following
method}\label{sc-Poly}
\setcounter{equation}{0}

Let us describe now our path-following predictor-corrector
strategy. It employs the following univariate function:
\beq\label{def-eta}
\ba{rcl}
\eta_{\bar \alpha}(\alpha) & = & \left\{ \ba{rl} 2
\alpha,& \alpha \in
[0, {1 \over 3} \bar \alpha),\\
{\alpha + \bar \alpha \over 2}, & \alpha \in [{1 \over 3}
\bar \alpha, \bar \alpha]. \ea \right.
\ea
\eeq
This function will be used for updating the length of the
current predictor step $\alpha$ by the rule $\alpha_+ =
\eta_{\bar \alpha}(\alpha)$. If the current step is small,
then it will be doubled. On the other hand, if $\alpha$ is
close enough to the maximal step size $\bar \alpha$, then
this distance for the new value $\alpha_+$ will be halved.

Let us fix the tolerance parameter $\beta'={1 \over 6}$
for the proximity measure $\Gamma$.
Consider the following method.
\beq\label{met-PC}
\ba{|c|}
\hline \\
\mbox{\bf Dual path-following method}\\
\mbox{\bf for
barriers with negative curvature}\\
\\
\ba{|rl|}
\hline & \\
{\bf 1.} & \mbox{Set $\mu_0=1$ and find point $y_0 \in
{\cal N}\left(\mu_0,{1 \over 25}
\right)$.}\\
& \\
{\bf 2.} & \mbox{\bf For $k \geq 0$ iterate:} \quad \mbox{a) Compute $\bar \alpha_k = \bar \alpha(y_k)$.}\\
& \\
& \mbox{b) Set $\alpha_{k,0} = {1 \over 6} \min \left\{1,
{1 \over  \|
v(y_k) \|_{y_k}} \right\}$. Using recurrence}\\
& \quad \alpha_{k,i+1} = \eta_{\bar
\alpha_k}(\alpha_{k,i}),
\mbox{ find the maximal $i \equiv i_k$,}\\
& \quad \mbox{such that
$\Gamma_{\mu_k}(y_k,\alpha_{k,i}) \leq \beta'$.}\\
& \\
& \mbox{c) Set $\alpha_k = \alpha_{k,i_k}$, $p_k = y_k +
\alpha_k v(y_k)$, $\mu_{k+1} =
{\mu_k \over \xi_{\bar \alpha_k}(\alpha_k)}$.}\\
& \\
&\mbox{d) Starting from $p_k$, apply the Newton method for}\\
&\quad \mbox{finding $y_{k+1} \in {\cal N}\left(\mu_{k+1},
\beta_{k+1}
\right)$ with $\beta_{k+1} = {\mu_{k+1} \over 25}$.}\\
& \\
\hline
\ea \\
\\
\hline
\ea
\eeq

Recall that the bound
$$
\ba{rcl}
\Gamma_{mu}(y,\alpha) & = & (1 + \alpha \sigma_{s(y)}(-A^*
v(y))) \| \nabla f(y(\alpha)) - {1 \over \mu} \xi_{\bar
\alpha}(\alpha) b \|_y
\ea
$$
is explicitly computable for different values of $\alpha$ (we
need to compute only the new vectors of the gradients
$\nabla f(y(\alpha))$).

Let us show now that the predictor-corrector scheme
(\ref{met-PC}) has polyno\-mi\-al-time complexity.
\BL\label{lm-Poly}
Suppose $F_*$ has negative curvature and let $y \in {\cal
N}(\mu,\beta)$ with $\beta \leq {1 \over 25}$. Then for
all
\beq\label{eq-Poly}
\ba{rcl}
\alpha \in \left[0, {1 \over 6} \min \left\{1, {1 \over \|
v(y_k) \|_{y_k}} \right\} \right]
\ea
\eeq
we have $\Gamma_{\mu}(y,\alpha) \leq \beta'$.
\EL
\proof
Denote $r = \| v(y) \|_y$, and $\hat r = \max\{ 1, \| v(y)
\|_y \}$. For every $\alpha \in [0,{1 \over 6 \hat r}]$ we
have
$$
\ba{l}
\Gamma_{\mu}(y,\alpha)\; \stackrel{(\ref{def-SH})}{\leq}
\; (1+\alpha r) \cdot \| \nabla
f(y(\alpha)) -
{\xi_{\bar \alpha}(\alpha) \over \mu} \cdot b \|_y\\
\\
\stackrel{(\ref{eq-Xi})}{\leq} \;
(1+\alpha r) \cdot
\left( \| \nabla f(y(\alpha)) - (1+\alpha) \nabla f(y)
\|_y + {\alpha^2 r \over \bar \alpha - \alpha }  +
\xi_{\bar \alpha}(\alpha) \| \nabla f(y) - {1 \over \mu}
\cdot b \|_y
\right)\\
\\
\stackrel{(\ref{eq-EstLoc}), (\ref{eq-DikLow})}{\leq}\;
 (1+\alpha r) \cdot
\left( {2\alpha^2 r^2 \over 1 - \alpha r}  + \beta \cdot
\left[
1 + {\alpha  \over 1 - \alpha r} \right] \right) \; \leq
\;  (1+\alpha \hat r) \cdot {2\alpha^2 \hat r^2 + \beta (1
+ \alpha) \over 1 - \alpha \hat r}.
\ea
$$
Hence, $\Gamma_{\mu}(y,\alpha)
\stackrel{(\ref{eq-Poly})}{\leq}
\left(1+ {1 \over 6} \right) \cdot
{ {2 \over 36 } + {1 \over 25}(1+{1 \over 6}) \,
 \over 1 - {1 \over 6}}
< {1 \over 6}.$
\qed

\BC
Let sequence $\left\{\mu_k\right\}$ be generated by method
(\ref{met-PC}). Then for every $k \geq 0$ we have
\beq\label{eq-MuRate}
\ba{rcl}
\mu_{k+1} & \leq &  \left(1 + {1 \over 6 \nu^{1/2}}
\right)^{-1} \mu_k.
\ea
\eeq
\EC
\proof
Note that by Lemma \ref{lm-Poly}, we get
$\Gamma_{\mu_k}(y_k, \alpha_{k,0}) \leq \beta'$.
Therefore, in method (\ref{met-PC}) we always have
$\alpha_k \geq \alpha_{k,0} \refGE{eq-NDec1} {1 \over 6
\nu^{1/2}}$. It remains to note that $\xi_{\bar
\alpha}(\alpha_k) \refGE{eq-Xi} 1 + \alpha_k$.
\qed

At the same time, it follows from Theorem \ref{th-NCG}
that method (\ref{met-PC}) can be accelerated. Taking into
account the choice $\beta_k = {\mu_k \over 25}$ in
(\ref{met-PC}), we see that
$$
\ba{rcl}
\kappa_1 & \leq & \nu + {{1 \over 25} ({1 \over 25} +
\nu^{1/2}) \over 1 - {1 \over 25}},\\
\\
\kappa_2 & \leq & {2 \over \gamma_d} (\sigma_d + {6 \over
25} \nu^2),
\ea
$$
and $\kappa = \kappa_1 \kappa_2$. Let us assume that
$\mu_k$ is small enough. Namely, we assume that
\beq\label{eq-2Bound}
\ba{rcl}
{1 \over 1 - \kappa \mu_k - 2 \beta_k} & \leq 2 \quad
\Leftrightarrow \quad \mu_k \; \leq \; {1 \over 2(\kappa +
{2 \over 25})}.
\ea
\eeq
Then
\beq\label{eq-Step3}
\ba{rcl}
\bar \alpha_k & \refLE{eq-BAlpha11} & 1 + {\kappa \mu_k \over
1 - \kappa \mu_k - {2 \over 25} \mu_k} \; \leq \; 1 +
2\kappa \mu_k \; \refLE{eq-2Bound} \; 2,\\
\\
\bar \alpha_k & \refGE{eq-BAlpha1} & {1 \over 1 + \kappa \mu_k
} \; \refGE{eq-2Bound} \; {2 \over 3}.
\ea
\eeq
Therefore, taking into account that
$$
\ba{rcl}
\sigma_{s(y)}(-A^* v(y)) & \refLE{def-SH} & \| A^* v(y)
\|_{s(y)} \; = \; \| v(y) \|_y \; \refLE{eq-NDec1} \;
\nu^{1/2},
\ea
$$
we have the following bound on the growth of the proximity
measure:
\beq\label{eq-GrowPM}
\ba{rcl}
\Gamma_{\mu_k}(y_k, \alpha) & \refLE{eq-Main} & (1 +
\bar \alpha_k \nu^{1/2})
\left[ {\alpha \bar \alpha_k \over (\bar \alpha_k -
\alpha)^2} \mu_k c_0 + {\mu_k \over 25} \left(1 + {\alpha
\bar \alpha_k \over \bar \alpha_k - \alpha} \right)
\right]\\
\\
& \refLE{eq-Step3} & (1 + 2 \nu^{1/2}) \mu_k
\left[ {\bar \alpha_k^2 \over (\bar \alpha_k -
\alpha)^2} c_0 + {1 \over 25} \left(1 + {2
\bar \alpha_k \over \bar \alpha_k - \alpha} \right)
\right],
\ea
\eeq
where $c_0 \refEQ{eq-2Bound} \kappa \nu^{1/2} + {2
\kappa_1 \over \gamma_d}
\left[ \sigma_d + {6 \over 25} \nu^2 + 2 \kappa \nu (1 +
2 \nu^{1/2}) {1 - {1 \over 25} \over 1 - {2 \over 25}}
\right]$.

Assume now that $\mu_k$ is even smaller, namely, that
\beq\label{eq-3Bound}
\ba{rcl}
\mu_k & \leq & {\beta' \over (1+2 \nu^{1/2}) ( 9c_0 + {7
\over 25})}.
\ea
\eeq
Denote by $\xi(\mu)$ the unique positive solution of the
equation
\beq\label{eq-Xi}
\ba{rcl}
c_0 \xi^2 + {1 \over 25} (1 + 2 \xi) & = & {\beta' \over
(1+2 \nu^{1/2}) \mu}.
\ea
\eeq
In view of assumption (\ref{eq-3Bound}), we have $1 \leq
{1 \over 3} \xi(\mu_k)$. Therefore,
\beq\label{eq-XiBound}
\ba{rcl}
{\beta' \over (1+2 \nu^{1/2}) \mu_k} & \leq & (c_0 +
{7\over 9 \cdot 25}) \xi^2(\mu_k).
\ea
\eeq
Note that for $\alpha(\mu_k)$ defined by the equation
$$
\ba{rcl}
{\bar \alpha_k \over \bar \alpha_k - \alpha(\mu_k)} & = &
\xi(\mu_k) \; \refGE{eq-3Bound} \; 3,
\ea
$$
we have $\Gamma_{\mu_k}(y_k,\alpha(\mu_k)) \leq \beta'$.
Therefore, $y_k + \alpha(\mu_k) v(y_k) \in Q$. At the same
time,
\beq\label{eq-BarAlph}
\ba{rcl}
\alpha(\mu_k) & \geq & {2 \over 3} \bar \alpha_k.
\ea
\eeq
Note that in view of the termination criterion of Step b)
in method (\ref{met-PC}) we either have $2 \alpha_k \geq
\alpha(\mu_k)$, or $\half (\bar \alpha_k + \alpha_k) \geq
\alpha(\mu_k)$. The first case is possible only if
$\alpha_k < {1 \over 3} \bar \alpha_k$. This implies
$\alpha(\mu_k) < {2 \over 3}\bar \alpha_k$, and this
contradicts the lower bound (\ref{eq-BarAlph}). Thus,
we have
$$
\ba{rcl}
\alpha_k & \geq & 2 \alpha(\mu_k) - \bar \alpha_k.
\ea
$$
Therefore,
$$
\ba{rcl}
\xi_{\bar \alpha_k}(\alpha_k) & = & 1 + {\alpha_k \bar
\alpha_k \over \bar \alpha_k - \alpha_k} \; \geq \; 1 +
{\bar \alpha_k (2 \alpha(\mu_k) - \bar \alpha_k) \over
2(\bar \alpha_k - \alpha(\mu_k))} \; \refGE{eq-BarAlph} \;
 1 +
{\bar \alpha_k^2 \over 6(\bar \alpha_k - \alpha(\mu_k))}\\
\\
& = & 1 + {1 \over 6} \bar \alpha_k \xi(\mu_k) \;
\refGE{eq-Step3} \; 1 + {1 \over 9} \xi(\mu_k)
\refGE{eq-XiBound} \; 1 + {1 \over 9}\sqrt{c_1 \over
\mu_k},
\ea
$$
where $c_1 = {\beta' \over (1 + 2 \nu^{1/2})(c_0+{7 \over
9 \cdot 25})}$.

Thus, we have proved the following theorem.
\BT
Suppose that $\mu_k$ satisfies condition
(\ref{eq-3Bound}). Then method (\ref{met-PC}) has a local
superlinear rate of convergence:
$$
\ba{rcl}
\mu_{k+1} & \leq  { 9  \over c_1^{1/2}} \; \mu_k^{3/2}.
\ea
$$
\ET

It remains to estimate the full complexity of one
iteration of method (\ref{met-PC}). It has two auxiliary
search procedures. The first one (that is Step b))
consists in finding an appropriate value of the predictor
step. We need an auxiliary statement on performance of its
recursive rule $\alpha_+ = \eta_{\bar \alpha}(\alpha)$.
\BL\label{lm-Pred}
If $\alpha \geq 0$ and $\alpha_+ = \eta_{\bar
\alpha}(\alpha)$, then $\xi_{\bar \alpha}(\alpha_+) \geq 2
\xi_{\bar \alpha}(\alpha) - 1$. Hence, for the recurrence
$$
\ba{rcl}
\alpha_{i+1} & = & \eta_{\bar \alpha}(\alpha_i), \quad i
\geq 0,
\ea
$$
we have $\xi_{\bar \alpha}(\alpha_i) \geq 1 + \alpha_0
\cdot 2^i$.
\EL
\proof
If $\alpha_+ = 2 \alpha$, then $\xi_{\bar
\alpha}(\alpha_+) = 1 + {2 \alpha \bar \alpha \over
\bar \alpha - 2 \alpha} \geq 1 + {2 \alpha \bar
\alpha \over \bar \alpha - \alpha} = 2\xi_{\bar
\alpha}(\alpha)-1$. If $\alpha_+ = {\alpha + \bar \alpha
\over 2}$, then
$$
\ba{rcl}
\xi_{\bar \alpha}(\alpha_+) & = & 1 + {\bar \alpha ( \bar
\alpha + \alpha ) \over
\bar \alpha - \alpha} \; = \; \xi_{\bar \alpha}(\alpha)
+ {\bar \alpha^2 \over
\bar \alpha - \alpha} \; \geq \; 2\xi_{\bar
\alpha}(\alpha)-1.
\ea
$$
Therefore, $\xi_{\bar \alpha}(\alpha_i) \geq 1 +
(\xi_{\bar \alpha}(\alpha_0)-1) \cdot 2^i
\stackrel{(\ref{eq-Xi})}{\geq} 1 + \alpha_0 \cdot 2^i$.
\qed

Note that the number of evaluations of the
proximity measure $\Gamma_{\mu_k}(y_k, \cdot)$ at Step b)
of method (\ref{met-PC}) is equal to $i_k+2$. Therefore,
for the first $N$ iterations of this method we have
$$
\ba{rcl}
\sum\limits_{k=0}^{N-1} (i_k+2) & \leq & 2N +
\sum\limits_{k=0}^{N-1} \log_2 {\mu_k \over \mu_{k+1}
\alpha_{k,0}} \; \leq \; N(2 + \half \log_2 \nu) - \log_2
\mu_N.
\ea
$$
Taking into account that for solving the problem
(\ref{prob-standD}) with absolute accuracy $\epsilon$ we
need to ensure $\mu_N \leq {\nu \over \epsilon}$, and
using the rate of convergence (\ref{eq-MuRate}), we
conclude that the total number of evaluations of the
proximity measure $\Gamma_{\mu_k}(y_k,\cdot)$
does not exceed $O(\nu^{1/2} \log_2 \nu \log_2{\nu \over \epsilon})$.

It remains to estimate the complexity of the correction
process (this is Step d)). This process cannot be too
long either. Note that the penalty parameters $\mu_k$ are bounded
from below by ${\epsilon \over \nu}$, where $\epsilon$ is
the desired accuracy of the solution. On the other hand,
the point $p_k$ belongs to the region of quadratic
convergence of the Newton method. Therefore, the number of
iterations at Step d) is bounded by $O(\ln \ln {\kappa_1
\over \epsilon})$. In Section \ref{sc-Disc}, we will
demonstrate on simple examples that the high accuracy in
approximating the trajectory of central path is crucial
for local superlinear convergence of the proposed
algorithm.

There exists another possibility for organizing the
correction process at Step d). We can apply a gradient
method in the metric defined by the Hessian of the
objective function at point $p_k$. Then the rate of
convergence will be linear, and we could have up to $O(\ln
{\nu \over \epsilon})$ correction iterations at each step
of method (\ref{met-PC}). However, each of such an
iteration will be cheap since it does not need evaluating
and inverting the Hessian.

\section{Discussion}\label{sc-Disc}
\setcounter{equation}{0}

\subsection{2D-examples}

Let us look now at several 2D-examples illustrating
different aspects of our approach. Let us start with the
following problem:
\beq\label{prob-Ex1}
\max\limits_{y \in \R^2} \{ \la b,y\ra: \; y_2 \geq 0, \;
y_1 \geq  y_2^2\}.
\eeq
For this problem, we can use the following barrier
function:
$$
\ba{rcl}
f(y) & = & - \ln (y_1 - y_2^2) - \ln y_2.
\ea
$$
We are going to check our conditions for the optimal point
$y_*=0$.

Problem (\ref{prob-Ex1}) can be seen as a restriction of
the following conic problem:
\beq\label{prob-Con1}
\ba{rl}
\max\limits_{s,y} \{ \la b, y \ra: & s_1 = y_1, \; s_2 =
y_2,\; s_3 = 1,\; s_4 = y_2,\quad s_1 s_3 \geq s_2^2,\quad
s_4 \geq 0 \},
\ea
\eeq
endowed with the barrier $F_*(s) = - \ln (s_1 s_3 - s_2^2)
- \ln s_4$. Note that
$$
\nabla F_*(s)  =
\left( {- s_3\over s_1 s_3 - s_2^2},
{2 s_2  \over s_1 s_3 - s_2^2}, {- s_1 \over s_1 s_3 -
s_2^2},  {-1 \over s_4} \right)^T.
$$
Denote $\omega =  s_1 s_3 - s_2^2$.  Since in problem
(\ref{prob-Con1}) $y_*=0$ corresponds to $s_* = e_3$, we
have the following representation:
$$
\nabla^2 F_*(s) s_* = \nabla^2 F_*(s) e_3 = { 1 \over
\omega^2} \cdot\left( s_2^2 , - 2 s_1 s_2, s_1^2, 0
\right)^T.
$$
Let us choose in the primal space the norm
$$
\ba{rcl}
\la B x, x \ra & = & x_1^2 + \half x_2^2 + x_3^2 + x_4^2.
\ea
$$
Then $\| \nabla^2 F_*(s) s_* \|_B = [s_1^2 + s_2^2 ]/
\omega^2$. Hence, the region $\| \nabla^2 F_*(s(y)) s_*
\|_B \leq \sigma_d$ is formed by vectors $y = (y_1,y_2)$
satisfying the inequality
$$
\ba{rcl}
y_1^2 + y_2^2 & \leq & \sigma_d (y_1 - y_2^2)^2.
\ea
$$
Thus, the boundary curve of this region is given by
equation
$$
\ba{rcl}
y_1 & = & y_2^2 + {1 \over \sigma_d^{1/2}} \sqrt{y_1^2 +
y_2^2},
\ea
$$
which has a positive slope $[\sigma_d-1]^{-\frac{1}{2}}$ at
the origin (see Figure 1). Note that the central path
corresponding to the vector $b = (-1,0)^T$ can be found
form the equations
$$
{1 \over y_1 - y_2^2} \; = \; {1 \over \mu}, \quad {1
\over y_2} \; = \; {2 y_2 \over y_1 - y_2^2}.
$$
Thus, its characteristic equation is $y_1 = 3 y_2^2$, and,
for any value of $\sigma_d$, it leaves the region of
quadratic convergence as $\mu \to 0$. It is interesting
that in our example, Assumption \ref{ass-Bound} is valid if
and only if the problem (\ref{prob-Ex1}) with $y^*=0$
satisfies Assumption \ref{ass-Sharp}.

\vspace{2ex}
\begin{minipage}[c]{12cm}
\setlength{\unitlength}{1.0cm}
\begin{picture}(10,10)
\put(1,0){\vector(0,1){10}} \put(1.2,9.5){$y_1$}
\put(0,1){\vector(1,0){10}} \put(9.5,0.5){$y_2$}

\put(0.5,0.5){$0$}

\linethickness{0.7mm} \qbezier(1,1)(5,1)(8.5,8)
\put(8.5,7.3){$y_1 = y_2^2$}
\put(1,1){\line(0,1){8.5}}

\linethickness{0.3mm} \qbezier(1,1)(4,1)(5,6) \thinlines

\put(2.5,6.3){Central path $b = (-1,0)^T$}

\qbezier(1,1)(6,3)(8,8.5)
\put(6,9){$y_1 = y_2^2 + {1 \over \sigma_d^{1/2}}
\sqrt{y_1^2 + y_2^2}$}
\put(1,1){\line(5,2){8}}
\put(7.5,4.5){$y_1 = y_2/\sqrt{\sigma_d-1}$}

\put(2.5,7.5){\fbox{$\| \nabla^2 F_*(s(y)) s_*
\|_B \leq \sigma_d$}}

\end{picture}

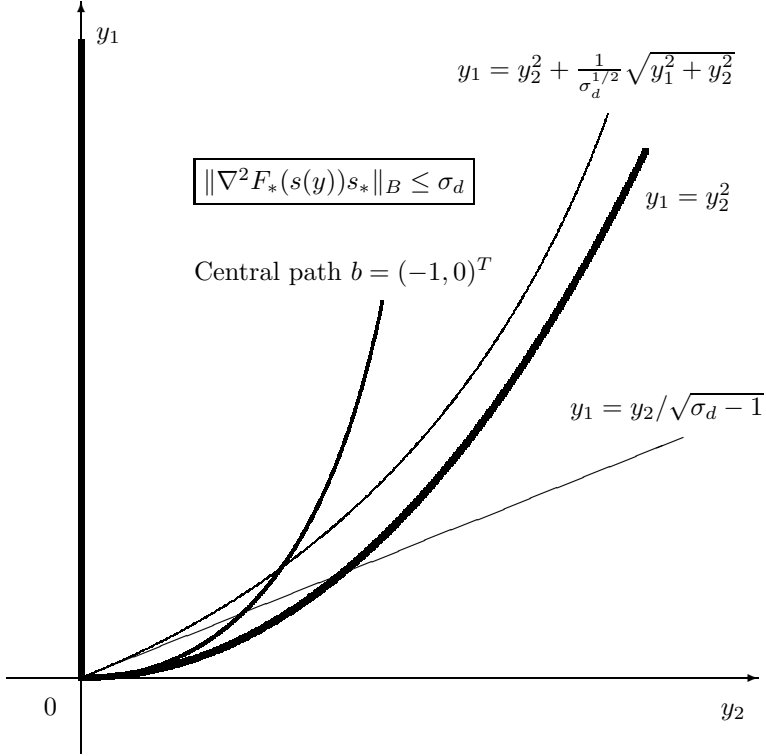
{\bf Figure 1.} Behavior of $\| \nabla^2 F_*(s(y))s_*
\|_B$.
\vspace{2ex}

\end{minipage}

In our second example we need the {\em maximal
neighborhood} of the central path:
\beq\label{def-Max}
\ba{rcl}
{\cal M}(\beta) & = & \cl \left( \bigcup\limits_{\mu \in
\R} {\cal N}(\mu,\beta)\right)  \\
\\
&  = & \left\{ y: \; \theta^2(y) \Def \| \nabla f(y)
\|_y^2 - {1 \over \| b \|_y^2} \la \nabla f(y), [\nabla^2
f(y) ]^{-1} b \ra^2 \leq \beta^2 \right\}.
\ea
\eeq
Note that $\theta(y) = \min\limits_{t \in \R} \| \nabla
f(y) - t b \|_y$.

Consider the following problem:
\beq\label{prob-Ex2}
\max\limits_{y \in \R^2} \{ y_1: \; \| y \| \leq 1 \}.
\eeq
where $\| \cdot \|$ is the standard Euclidean norm. Let us
endow the feasible set of this problem with the standard
barrier function $f(y) = - \ln (1 - \| y \|^2)$. Note that
$$
\ba{rcl}
\nabla f(y) & = & {2 y \over 1 - \| y \|^2}, \quad
\nabla^2 f(y) \; = \; {2 I \over 1 - \| y \|^2} + {4 y y^T
\over (1 - \| y \|^2)^2},\\
\\
\left[ \nabla^2 f(y) \right]^{-1} & = &
{1 - \| y \|^2 \over 2} \left( I - {2 y y^T \over 1 + \| y
\|^2} \right),\quad
\left[ \nabla^2 f(y) \right]^{-1} \nabla f(y) \; = \; {1 -
\| y \|^2 \over 1 + \| y \|^2} \cdot y.
\ea
$$
Therefore,
$$
\ba{c}
\| \nabla f(y) \|_y^2 \; = \; {2 \| y \|^2 \over 1 +\| y
\|^2},
\ea
$$
and for $b=(1,0)^T$ we have
$$
\ba{c}
\| b \|_y^2 \; = \; {1 - \| y \|^2 \over 2}
\cdot {1 - y_1^2 + y_2^2 \over 1 + \| y \|^2},
\quad
\la \nabla f(y),  \left[ \nabla^2 f(y) \right]^{-1} b \ra
\; = \; {1 - \| y \|^2 \over 1 + \| y \|^2} \cdot y_1.
\ea
$$
Thus,
$$
\ba{rcl}
\theta^2(y) & = & {2 \| y \|^2 \over 1 + \| y \|^2} - {2
\over 1 - \| y \|^2} \cdot {1 + \| y \|^2 \over 1 - y_1^2
+ y_2^2} \cdot  {(1 - \| y \|^2)^2 y_1^2 \over (1 + \| y
\|^2)^2} \\
\\
& = & {2  \over 1 + \| y \|^2} \left( \| y \|^2 - {y_1^2(1
- \| y \|^2)  \over 1 - y_1^2 + y_2^2} \right) \; = \; {2
y_2^2 \over 1 - y_1^2 + y_2^2}.
\ea
$$
We conclude that for problem (\ref{prob-Ex2}) the maximal
neighborhood of the central path has the following
representation:
\beq\label{eq-MaxRep}
\ba{rcl}
{\cal M}(\beta) & = & \left\{ y \in \R^2: \; y_1^2 + {2 -
\beta^2 \over \beta^2} \cdot y_2^2 \; \leq \; 1 \right\}
\ea
\eeq
(see Figure 2).

\vspace{2ex}
\begin{minipage}[c]{12cm}
\setlength{\unitlength}{1.0cm}
\begin{picture}(10,10)
\put(5,0){\vector(0,1){10}} \put(5.2,9.5){$y_1$}
\put(0,5){\vector(1,0){10}}
\put(9.5,5.2){$y_2$}
\put(5,5){\circle{8}}
\qbezier(4,5)(4,9)(5,9) \qbezier(5,9)(6,9)(6,5)
\qbezier(6,5)(6,1)(5,1) \qbezier(5,1)(4,1)(4,5)

\qbezier(4.6,5)(4.6,9)(5,9) \qbezier(5,9)(5.4,9)(5.4,5)
\qbezier(5.4,5)(5.4,1)(5,1) \qbezier(5,1)(4.6,1)(4.6,5)

\linethickness{1mm}
\put(5,1){\line(0,1){8}}
\thinlines
\put(5.1,4.6){$0$}
\put(6.3,4){\fbox{Central path}}
\put(6.3,4.1){\vector(-1,0){1.3}}
\put(6.3,0.4){\fbox{Large neighborhood}}
\put(6.3,0.78){\vector(-1,1){0.83}}
\put(0.44,0.4){\fbox{Small neighborhood}}
\put(4.03,0.78){\vector(1,1){0.83}}
\put(6.3,9){{\bf Problem:} $\quad\max\limits_{\| y \| \leq 1} y_1$.}
\put(5,5){\line(-1,6){0.63}}
\put(4.63,7.24){\circle*{0.1}}
\put(4.3,7){$y$}
\put(4.41,8.5){\circle*{0.1}}
\put(4.76,8.5){\circle*{0.1}}
\put(3.8,8.2){$p'$}
\put(4.41,8.5){\vector(1,0){0.37}}
\put(4.36,8.75){\circle*{0.1}}
\put(4,9.2){$p(y)$}
\put(5.8,8.4){$y_+$}
\end{picture}

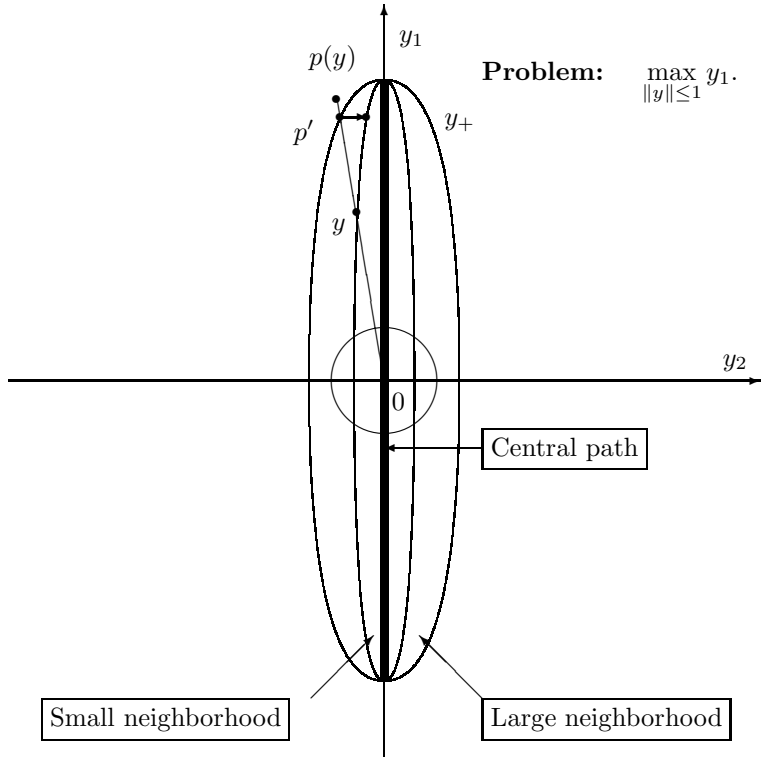
{\bf Figure 2.} Prediction in the absence of sharp
maximum.
\vspace{2ex}
\end{minipage}

Note that $p(y) = {2 y \over 1 + \| y \|^2} \in \inter Q$.
If the radii of the small and large neighborhoods of the
central path are fixed, by straightforward computations we
can see that the simple predictor-corrector update $y \to
y_+$ shown in Figure 2 has local linear rate of
convergence. In order to get a superlinear rate, we need to
tighten the small neighborhood of the central path as $\mu
\to 0$.

\subsection{Examples of cones with negative curvature}

In accordance with the definition (\ref{eq-NC}), negative
curvature of barrier functions is preserved by the
following operations.
\begin{itemize}
\item
If barriers $F_i$ for cones $K_i \subset \E_i$, $i = 1,2$, have
negative curvature, then the curvature of the barrier
$F_1+F_2$ for the cone $K_1 \oplus K_2$ is negative.
\item
If barriers $F_i$ for cones $K_i \subset \E$, $i = 1,2$, have
negative curvature, then the curvature of the barrier
$F_1+F_2$ for the cone $K_1 \bigcap K_2$ is negative.
\item
If barrier $F$ for cone $K$ has negative curvature, then
the curvature of the barrier $f(y) = F(A^*y)$ for the cone
$K_y = \{ y \in \H: \; A^* y \in K \}$ is negative.
\item
If barrier $F(x)$ for cone $K$ has negative curvature,
then the curvature of its restriction onto the linear
subspace $\{ x \in \E: \; A x = 0 \}$ is negative.
\end{itemize}
At the same time, we know two important families of cones
with negative curvature.
\begin{itemize}
\item
Self-scaled barriers have negative curvature (see
Corollary 3.2(i) in \cite{NT1}).
\item
Let $p(x)$ be hyperbolic polynomial. Then the barrier
$F(x) = - \ln p(x)$ has negative curvature (see \cite{Guler1997}).
\end{itemize}
Thus, using above mentioned operations, we can construct
barriers with negative curvature for many interesting
cones. In some situations, we can argue that currently,
some nonsymmetric
treatments of the primal-dual problem pair have better complexity
bounds than the primal-dual symmetric treatments.
\BE
Consider the cone of nonnegative polynomials:
$$
\ba{rcl}
K & = & \left\{ p \in \R^{2n+1}: \; \sum\limits_{i=0}^{2n}
p_i t^i \geq 0, \; \forall t \in \R \right\}.
\ea
$$
The dual to this cone is the cone of positive semidefinite
Hankel matrices. For $k \in \{0,1, \dots , 2n\}$, denote
$$
\ba{rcl}
H_k \in \R^{(n+1)\times (n+1)} : \; H_k^{(i,j)} & = &
\left\{ \ba{rl} 1, &\mbox{if $i+j=k+2$} \\ 0,
&\mbox{otherwise} \ea \right., \quad i,j \in \{0,1, \dots, n\}.
\ea
$$
For $s \in \R^{2n+1}$ we can define now the following
linear operator:
$$
\ba{rcl}
H(s) & = & \sum\limits_{i=0}^{2n} s_i \cdot H_i.
\ea
$$
Then the cone dual to $K$ can be represented as follows:
$$
\ba{rcl}
K^* & = & \{ s \in \R^{2n+1}:\; H(s) \succeq 0 \}.
\ea
$$

The natural barrier for the dual cone is $f(s) = - \ln
\det H(s)$. Clearly, it has negative curvature. Note that
we can lift the primal cone to a higher dimensional space
(see \cite{PP}):
$$
\ba{rcl}
K & = & \left\{ p \in \R^{2n+1}: \; p_i = \la H_i, Y \ra,\; Y
\succeq 0,\; i \in \{0,1,\dots, 2n\} \right\},
\ea
$$
and use $F(Y) = -\ln \det Y$ as a barrier function for the
extended feasible set. However, in this case we
significantly increase the number of variables. Moreover,
we need $O(n^3)$ operations for computing the value of the
barrier $F(Y)$ and its gradient. On the other hand, in the
dual space the cost of all necessary computations is very
low ($O(n \ln^2 n)$ for the function value and $O(n^2
\ln^2 n)$ for solution of the Newton system, see
\cite{PP-Comp}). On top of these advantages, for
non-degenerate dual problems, now we have a locally
superlinearly convergent path-following scheme
(\ref{met-PC}).
\EE

To conclude the paper, let us mention that the negative
curvature seems to be a natural property of
self-concordant barriers. Indeed, let us move from some
point $x \in \inter K$ along the direction $h \in K$: $u = x
+h$. Then the Dikin ellipsoid of barrier $F$ at point $x$,
moved to the new center $u$, still belongs to $K$:
$$
\ba{rcl}
u + (W_r(x)-x) & = & h + W_r(x) \; \subset \; K.
\ea
$$
We should expect that, in this situation, the Dikin ellipsoid $W_r(u)$
becomes even larger (in any case, we should expect that it does
not get smaller). This is exactly the negative curvature
condition: $\nabla^2 F(x) \succeq \nabla^2 F(u)$.
Thus, we are led to the following unsolved problem:
\emph{What is the class of regular convex cones that admit a logarithmically
homogeneous self-concordant barrier with negative curvature?
In particular, does every regular convex cone admit a logarithmically
homogeneous self-concordant barrier with negative curvature?}


\begin{thebibliography}{NN}

\bibitem{Bauschke2001}
H. H. Bauschke, O. G\"{u}ler, A. S. Lewis and H. Sendov,
Hyperbolic polynomials and convex analysis,
\emph{Canad. J. Math.} 53 (2001) 470--488.

\bibitem{PP-Comp}
Y.~Genin, Y.~Hachez, Yu.~Nesterov and P.~Van Dooren,
Optimization problems over positive pseudo-polynomial
matrices, \emph{SIAM J. Matrix Anal.
Appl.} 25 (2003) 57--79.

\bibitem{FM}
A. V. Fiacco and G. P. McCormick, \emph{Nonlinear Programming:
Sequential Unconstrained Minimization Techniques}, J.Wiley
$\&$ Sons Inc., New York, 1968.


\bibitem{Guler1997}
O. G\"{u}ler,
Hyperbolic polynomials and interior-point methods
for convex programming,
\emph{Math. Oper. Res.} 22 (1997) 350--377.


\bibitem{JPS1999}
J. Ji, F. A. Potra and R. Sheng,
On the local convergence of a predictor-corrector
method for semidefinite programming,
\emph{SIAM J. Optim.} 10 (1999) 195--210.

\bibitem{KSS1998}
M. Kojima, M. Shida and S. Shindoh,
Local convergence of predictor-corrector
infeasible-interior-point algorithms for SDPs and SDLCPs,
\emph{Math. Program.} 80 (1998) 129--160.

\bibitem{LuMo2007}
Z. Lu and R. D. C. Monteiro,
Limiting behavior of the Alizadeh-Haeberly-Overton
weighted paths in semidefinite programming,
\emph{Optim. Methods Softw.} 22 (2007) 849--870.


\bibitem{LuMo2004}
Z. Lu and R. D. C. Monteiro,
Error bounds and limiting behavior of
weighted paths associated with
the SDP map $X^{1/2}SX^{1/2}$,
\emph{SIAM J. Optim.} 15 (2004/05) 348--374.

\bibitem{LSZ1998}
Z.-Q. Luo, J. F. Sturm and S. Zhang,
Superlinear convergence of a symmetric
primal-dual path following algorithm for
semidefinite programming,
\emph{SIAM J. Optim.} 8 (1998) 59--81.

\bibitem{Meh1993}
S. Mehrotra,
Quadratic convergence in a primal-dual method,
\emph{Math. Oper. Res.} 18 (1993) 741--751.

\bibitem{Intro}
Yu. Nesterov, {\em Introductory Lectures on Convex
Optimization}, Kluwer, Boston, 2004.

\bibitem{PP}
Yu. Nesterov, Squared functional systems and optimization
problems, In {\sl High Performance Optimization}, H.Frenk,
T.Terlaky and S.Zhang (Eds.), Kluwer, 1999. pp.405--439.

\bibitem{NN}
Yu. Nesterov and A. Nemirovskii, {\em Interior Point
Polynomial Methods in Convex Programming: Theory and
Applications,\/} SIAM, Philadelphia, 1994.

\bibitem{NT1}
Yu. Nesterov and M. J. Todd, Self-scaled barriers and
interior-point methods for convex programming,
\emph{Math. Oper. Res.} 22 (1997) 1--42.

\bibitem{PS2000}
F. A. Potra and R. Sheng,
Superlinear convergence of a
predictor-corrector method for
semidefinite programming without
shrinking central path neighborhood,
\emph{Bull. Math. Soc. Sci. Math. Roumanie}
43 (2000) 107--124.

\bibitem{PS1998}
F. A. Potra and R. Sheng,
Superlinear convergence of interior-point
algorithms for semidefinite programming,
\emph{J. Optim. Theory Appl.} 99 (1998) 103--119.

\bibitem{Renegar2006}
J. Renegar,
Hyperbolic programs, and their derivative relaxations,
\emph{Found. Comput. Math.} 6 (2006) 59--79.

\bibitem{YGTZ1993}
Y. Ye, O. G\"{u}ler, R. A. Tapia and Y. Zhang,
A quadratically convergent $O(\sqrt{n}L)$-iteration
algorithm for linear programming,
\emph{Math. Program.} 59 (1993) 151--162.

\bibitem{ZTD1992}
Y. Zhang, R. A. Tapia and J. E. Dennis,
On the superlinear and quadratic convergence of
primal-dual interior point linear programming algorithms,
\emph{SIAM J. Optim.} 2 (1992) 304--324.



\end{thebibliography}
\end{document}